\documentclass[11pt]{article}
\usepackage{pxfonts}
\usepackage{yfonts}
\usepackage{dsfont}
\usepackage{graphicx}
\usepackage{relsize}

\def\bel{\begin{equation}\label}
\def\eeq{\end{equation}}
\def\ds{\displaystyle}
\def\endproof{\hphantom{MM}
\hfill\llap{$\square$}\goodbreak}

\def\mt{\longrightarrow}
\def\v{\vskip 1em}

\def\R{\mathbb R}
\def\Z{\mathbb Z}
\def\C{\mathfrak{C}}

\def\Re{{\bf Re}}
\def\Im {{\bf Im}}
\def\S{{\bf S}}

\def\F{\mathfrak{F}}

\def\O{{\bf O}}
\def\Q{\mathfrak{Q}}

\def\B{{\bf B}}
\def\H{{\bf H}}
\def\L{{\bf L}}

\def\BMO{{\bf BMO}}

\def\p{{\partial}}

\def\i{{\bf i}}

\def\Hat{\widehat}

\def\bar{\overline}
\def\supp{{\bf supp}}

\def\M{{\bf M}}

\def\alpha{\alphaup}
\def\beta{\betaup}
\def\gamma{\gammaup}
\def\delta{\deltaup}
\def\xi{{\xiup}}
\def\eta{{\etaup}}
\def\tau{{\tauup}}
\def\rho{{\rhoup}}
\def\phi{{\phiup}}
\def\psi{{\psiup}}
\def\lambda{{\lambdaup}}
\def\omega{\omegaup}
\def\varphi{{\varphiup}}
\def\gamma{{\gammaup}}
\def\c{{\bf c}}

\newtheorem{thm}{Theorem}[section]

\newtheorem{remark}{Remark}[section]

\begin{document}
 \[\begin{array}{cc}\hbox{\LARGE{\bf Regularity of Fourier integral operators on product spaces}}
 \end{array}\]

 \[\hbox{Zipeng Wang}\]
 \[\begin{array}{ccc}
 \hbox{Department of Mathematics, Westlake university}
 \\ 
 \hbox{Cloud town, Hangzhou of China}
 \end{array}\]
 \begin{abstract}
We study the regularity of Fourier integral operators, by allowing their symbols  to  satisfy certain multi-parameter characteristics. As a result, we prove a sharp $\L^p$-estimate obtained by  Seeger, Sogge and Stein  on product spaces.
\end{abstract}

\section{Introduction}
\setcounter{equation}{0}
 Let $f$ be a Schwartz function.
 We consider a Fourier integral operator $\F$  defined by
\bel{Ff}
\Big(\F f\Big)(x)~=~\int_{\R^n} f(y)\Omega(x,y)dy
\eeq
whose kernel is given by an oscillatory integral
\bel{Kernel}
\Omega(x,y)~=~\int_{\R^n}e^{2\pi\i \left(\Phi(x,\xi)-y\cdot\xi\right)}\sigma(x,y,\xi)d\xi.
\eeq
We require that the symbol  $\sigma(x,y,\xi)\in\mathcal{C}^\infty(\R^n\times\R^n\times\R^n)$ has a compact support in  $x$ and $y$. 
 On the other hand,  the phase function $\Phi(x,\xi)$ is real, homogeneous of degree $1$ in $\xi$, smooth for every $x$ and $\xi\neq0$. Moreover, it satisfies   
 the nondegeneracy condition
\bel{nondegeneracy}
 \det\left({\p^2\Phi\over \p x_i\p \xi_j}\right)\left(x,\xi\right)~\neq~0,\qquad \xi~\neq~0
 \eeq
on the support of $\sigma(x,y,\xi)$.  

Let $\C$ denote a generic constant with subindices indicating its dependence.   
 We say $\sigma\in \hbox{S}^m$ if 
\bel{class}
\left|\p_\xi^\alpha\p_{x,y}^\beta \sigma(x,y,\xi)\right|~\leq~\C_{\alpha~\beta}\left(1+|\xi|\right)^m \left({1\over 1+|\xi|}\right)^{|\alpha|}
\eeq
for every multi-indices $\alpha,\beta$.
 
Fourier integral operator $\F$  defined in (\ref{Ff})-(\ref{nondegeneracy}) has been extensively studied since the 1970's for its own right of interest found in harmonic analysis. 
Let $\sigma\in\hbox{S}^0$. The $\L^2$-boundedness of $\F$  and also its generalization associated to some appropriate local canonical graph,  was shown by Eskin \cite{Eskin} and H\"{o}rmander \cite{Hormander}. In contrast to this $\L^2$-estimate, it is well known that Fourier integral operator of order zero is not bounded on $\L^p$-spaces for $p\neq2$. The optimal $\L^p$-result was first investigated by H\"{o}rmander \cite{Hormander}, then Duistermaat and H\"{o}rmander \cite{Duistermaat-Hormander}
and eventually proved by Seeger, Sogge and Stein \cite{S.S.S}.

{\bf Theorem A: ~Seeger, Sogge and Stein, 1991}\\
 {\it Let} $\sigma\in \hbox{S}^m$ {\it as (\ref{class}) for $-(n-1)/2<m\leq0$. Fourier integral operator $\F$   defined in (\ref{Ff})-(\ref{nondegeneracy}) extends to a bounded operator 
  \bel{result}\left\| \F f\right\|_{\L^p(\R^n)}~\leq~\C_{p~\sigma~\Phi}~\left\| f\right\|_{\L^p(\R^n)}
 \eeq
whenever 
\bel{formula}
\left| {1\over 2}-{1\over p}\right|~\leq~{-m\over n-1}.
\eeq}
\begin{remark} {\bf Theorem A} is sharp: Let $a(x)b(y)\in\mathcal{C}^\infty_o(\R^n\times\R^n)$ where $a(x)\neq0$ for $|x|=1$ and $b(y)\equiv1$ for $|y|<1$.  Observe that $\sigma(x,y,\xi)\doteq a(x)b(y)\left(1+|\xi|\right)^m\in\hbox{S}^m$. Suppose $\Phi(x,\xi)=x\cdot\xi+|\xi|$. Fourier integral operator $\F$ is not bounded on $\L^p(\R^n)$ for $\left|1/2-1/p\right|>-m/(n-1),~(1-n)/2\leq m\leq0$.
\end{remark}
A regarding estimate of {\bf Remark 1.1} is given by {\bf 6.13}, chapter IX of Stein \cite{Stein} whereas the result is obtained by using the asymptotic of Bessel functions. Note that  Fourier integral operator with  phase function $\Phi(x,\xi)=x\cdot\xi\pm|\xi|$ arose to solve the wave equation, as was investigated by Colin de Verdi\'{e}re and Frisch \cite{Colin-Frisch}, Beals \cite{Beals}, Brenner \cite{Brenner}, Peral \cite{Peral} and Miyachi \cite{Miyachi}.

In this paper, we give an extension of the sharp $\L^p$-estimate in {\bf Theorem A}, by studying Fourier integral operator whose symbol satisfies a 2-parameter characteristic. 

Let $\xi=(\tau,\lambda)\in\R\times\R^{n-1}$.
We say $\sigma\in\S^m$ if 
\bel{Class}
\left|\p_\tau^\alpha\p_{\lambda}^\beta \p_{x,y}^\gamma\sigma(x,y,\xi)\right|~\leq~\C_{\alpha~\beta~\gamma}\left(1+|\xi|\right)^m\left({1\over 1+|\tau|}\right)^{\alpha}\left({1\over 1+|\lambda|}\right)^{|\beta|}
\eeq
for every multi-indices $\alpha,\beta$ and $\gamma$. 

Study of certain operators  that commute with a multi-parameter family of dilations  dates back to the time of Jessen, Marcinkiewicz and Zygmund.  Over the past several decades,
 a number of pioneering  results  have been accomplished,  for example 
by Robert Fefferman \cite{R.Fefferman}-\cite{R.Fefferman''}, Chang and  Fefferman \cite{Chang-Fefferman}, Cordoba and Fefferman \cite{Cordoba-Fefferman}, Fefferman and Stein \cite{R-F.S},
M\"{u}ller, Ricci and Stein \cite{M.R.S},
Journ\'{e} \cite{Journe'} and Pipher \cite{Pipher}. Our main result is the following.
\v

{\bf Theorem A*}~~
 {\it Let} $\sigma\in \S^m$ {\it as (\ref{Class}) for $-(n-1)/2<m\leq0$. Fourier integral operator $\F$   defined in (\ref{Ff})-(\ref{nondegeneracy}) extends to a bounded operator 
  \bel{result}\left\| \F f\right\|_{\L^p(\R^n)}~\leq~\C_{p~\sigma~\Phi}~\left\| f\right\|_{\L^p(\R^n)}
 \eeq
whenever 
\bel{Formula*}
\left| {1\over 2}-{1\over p}\right|~\leq~{-m\over n-1}.
\eeq}

{\bf Theorem A*} can be further generalized to an $n$-parameter setting. We leave the discussion to the end of this paper.
Some more recent works   in the direction of  harmonic analysis  on product spaces refer to Tanaka and Yabuta \cite{Tanaka-Yabuta}, Sawyer and Wang \cite{Sawyer-Wang}-\cite{Sawyer-Wang'} and Wang \cite{Wang}. 
 Historical background of Fourier integral operators  can be found in the books by Sogge \cite{Sogge} and Duistermaat \cite{Duistermaat}.

\section{Cone decomposition}
\setcounter{equation}{0}
We introduce a new framework where the frequency space is decomposed into an infinitely many dyadic cones. 
Every partial  operator whose symbol is supported on a dyadic cone   essentially  is a one-parameter Fourier integral operator, satisfying the desired regularity. Moreover, its norm decays exponentially as the eccentricity  of the cone getting large.

Let $\varphi$ be a smooth {\it bump}-function such that 
\bel{varphi funda}
\varphi(\xi)~\equiv1~~~~\hbox{for} ~~~~|\xi|~\leq~1\qquad\hbox{and} \qquad\varphi(\xi)~=~0~~~~\hbox{for}~~~~ |\xi|>2.
\eeq
Consider
\bel{delta_t}
\deltaup_\ell(\xi)~=~\varphi\left(2^{\ell}{\tau\over|\lambda|}\right)-\varphi\left(2^{\ell+1}{\tau\over|\lambda|}\right),\qquad \ell\in\Z
\eeq
which is supported on the dyadic cone
\bel{Cone}
\Lambda_\ell~=~ \left\{  (\tau,\lambda)\in\R\times\R^{n-1}~\colon~      2^{-\ell-1}~<~ {|\tau|\over|\lambda|}~<~2^{-\ell+1} \right\}.
\eeq

\begin{figure}[h]
\centering
\includegraphics[scale=0.40]{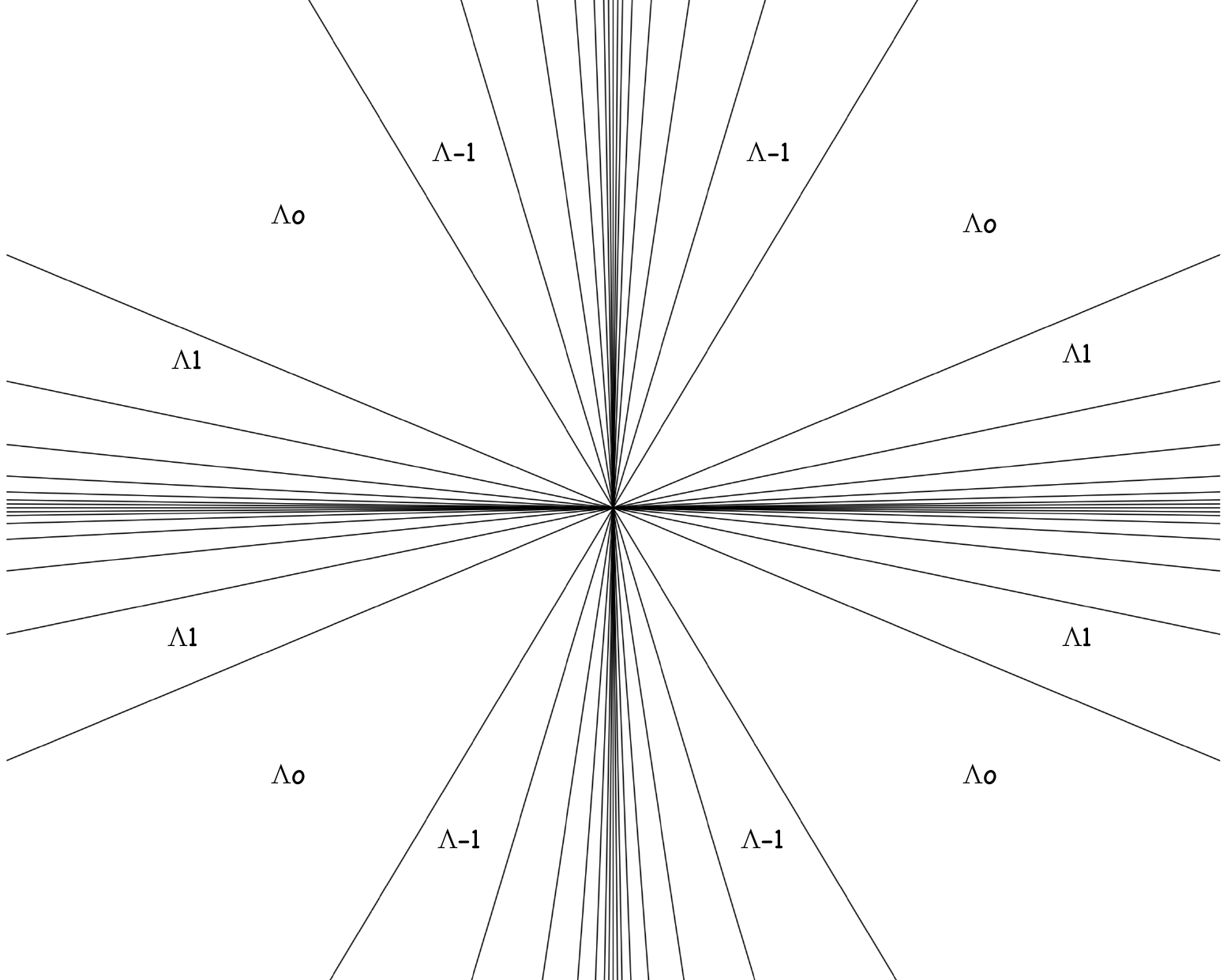}
\caption{$\tau$ is in the vertical direction and $\lambda$ is in the horizontal direction.}
\end{figure}
 
Define
\bel{Partial}
\begin{array}{cc}\ds
\Big(\F_\ell f\Big)(x)~=~\int_{\R^n} f(y)\Omega_\ell(x,y)dy,
\\\\ \ds
\Omega_\ell(x,y)~=~\int_{\R^n}e^{2\pi\i \left(\Phi(x,\xi)-y\cdot\xi\right)}\sigma(x,y,\xi)\delta_\ell(\xi)d\xi.
\end{array}
\eeq
Recall that $\Phi(x,\xi)$ is real, homogeneous of degree $1$ in $\xi$ and smooth for every $x$ and $\xi\neq0$, satisfying the nondegeneracy condition in (\ref{nondegeneracy}). 
\v
{\bf Theorem A**} ~{\it Let $\sigma\in \S^m$  as (\ref{Class}) for $-(n-1)/2<m\leq0$. Every $\F_\ell$ defined in (\ref{delta_t})-(\ref{Partial}) extends to a bounded operator 
\bel{Decay Result>}
\left\| \F_\ell f\right\|_{\L^p(\R^n)}~\leq~\C_{p~\sigma~\Phi}~2^{\left({m\over n}\right)\ell} \left\| f\right\|_{\L^p(\R^n)},\qquad \ell\ge0
\eeq
and
\bel{Decay Result<}
  \left\| \F_\ell f\right\|_{\L^p(\R^n)}~\leq~\C_{p~\sigma~\Phi}~2^{-\left({m\over n}\right)(n-1)\ell} \left\| f\right\|_{\L^p(\R^n)},\qquad \ell\leq0 
\eeq
whenever
\bel{Formula**}
\left| {1\over 2}-{1\over p}\right|~\leq~{-m\over n-1}.
\eeq}
\v
It is clear that {\bf Theorem A**} implies {\bf Theorem A*} for $p\neq2$ by applying Minkowski inequality.

{\bf Sketch of Proof:} First, in section 3, we prove  $\F$  ( also $\F_\ell$ ) bounded on $\L^2(\R^n)$  and show $\F_\ell\colon\L^p(\R^n)\mt\L^2(\R^n)$ for $-m/n=1/p-1/2$ and $\F_\ell\colon\L^2(\R^n)\mt\L^{q}(\R^n)$ for $-m/n=1/2-1/q$ with desired operator norms. Then, at every  partial operator $\F_\ell$, we  develop our analysis  in the same spirt of Seeger, Sogge and Stein \cite{S.S.S}. In particular, we   study its regularity  on   $\H^1$-Hardy space, by considering  separately for $\F_\ell$ restricted to the so-called  {\it region of influence}, denoted by  $\Q_r$, and its complement subset ${^c}\Q_r=\R^n\setminus\Q_r$. 

In section 4, by using the $\L^2$-estimates obtained in section 3, we prove $\F_\ell\colon\H^1(\R^n)\mt\L^1(\Q_r)$ for $\sigma\in\S^{-{n-1\over 2}}$   and then  give an heuristic estimate of (\ref{Decay Result>})-(\ref{Decay Result<}) by an interpolation argument.
In section 5, we show  $\F_\ell\colon\H^1(\R^n)\mt\L^1({^c}\Q_r)$  provided that the kernel of $\F_\ell$ in (\ref{Partial}) satisfies certain majorization properties, accumulated in the lemma named as {\bf Principal Lemma}.

From section 6, we begin to construct a second dyadic decomposition where the frequency space is decomposed into finitely many geometric cones $\Gamma^\nu_j$ for $j\ge0$ fixed whose central direction $\xi^\nu_j$   are almost equally distributed on the unit sphere $\mathbb{S}^{n-1}$ (as $\nu$ varies).  Their intersection with the dyadic annuli $\{2^{j-1}\leq|\xi|\leq2^{j+1}\}$ form a collection of thin rectangles. The corresponding  partial operators  having their  symbols  supported on these  thin rectangles are essentially the  non-isotropic Calder\'{o}n-Zygmund operators. Furthermore,  their norms can be added up for $\sigma \in\S^{-{n-1\over 2}}$.  A similar idea was used by Fefferman \cite{FC1},  C\'{o}rdoba \cite{Cordoba} and  Christ and Sogge \cite{Christ and Sogge}  in study of Bochner-Riesz multipliers.

On the other hand, note that $\F_\ell$ itself defined in (\ref{delta_t})-(\ref{Partial}) has a symbol whose support is restricted to the dyadic cone $\Lambda_\ell$  in (\ref{Cone}). The crucial part of our analysis is to study those partial operators having symbols  supported on the intersection  $\Gamma^\nu_j\cap\Lambda_\ell\cap\{2^{j-1}\leq|\xi|\leq2^{j+1}\}$. For $j\ge0$ fixed and $\ell$  changing, situations must be handled in different ways.
For this reason, we give a $3$-fold estimate with respect to
 $-j/2-3\leq\ell\leq j/2+2$, $\ell>j/2+3$ and $\ell<-j/2-3$,
and prove {\bf Principal Lemma} in section 7.
 
\begin{remark}
The proof will be self-contained, except for 
(\ref{Phi_x est}) and (\ref{d Est Psi}),  which are given explicitly  in {\it 3.1.1}, chapter IX and {\bf 4.5}, chapter IX of the book by Stein \cite{Stein}.

We like to emphasize that both (\ref{Phi_x est}) and (\ref{d Est Psi}) are only associated to the phase function $\Phi(x,\xi)$ satisfying those essential conditions: real, homogeneous of degree $1$ in $\xi$,   smooth for every $x$ and $\xi\neq0$, and the nondegeneracy condition in (\ref{nondegeneracy})    but has no further restriction added in this paper.
\end{remark}

\section{$\L^2$-boundedness of Fourier integral operators}
\setcounter{equation}{0}
We first show that $\F$ defined in (\ref{Ff})-(\ref{Kernel}) can be written as a finite sum of 
\bel{Ff'}
\Big(\mathcal{F} f\Big)(x)~=~\int_{\R^n}e^{2\pi\i\Phi(x,\xi)}\sigma(x,\xi)\Hat{f}(\xi)d\xi
\eeq
where $\Hat{f}$ is the Fourier transform of $f$ and $\sigma(x,\xi)\in\mathcal{C}^\infty(\R^n\times\R^n)$ has a compact support in $x$.

Recall that $\sigma(x,y,\xi)\in\mathcal{C}^\infty(\R^n\times\R^n\times\R^n)$  has a compact support in $x$ and $y$. On the $y$-space, we can construct a smooth partition of unity $\left\{\phi^\nu\right\}_\nu$ such that every $\phi^\nu$ is supported on a ball $B_r(y^\nu)$ centered on $y^\nu\in\R^n$ with radius $r\leq{1/ 2n}$. Moreover, there are only finitely many $\nu$ s, depending on the size of $y$-support of $\sigma(x,y,\xi)$.

For each $\sigma(x,y,\xi)\phi^\nu(y)$, we write its Taylor expansion $w.r.t~y$ centered on $y^\nu$.
Let $\hbox{\bf 1}_{\nu}(y)$ denote the indicator function on $B_r(y^\nu)$. Observe that
\bel{Ff Sum Taylor}
\begin{array}{lr}\ds
\int_{\R^n} f(y)\left\{\int_{\R^n}e^{2\pi\i \left(\Phi(x,\xi)-y\cdot\xi\right)}\sigma(x,y,\xi)\phi^\nu(y)d\xi\right\}dy
\\\\ \ds
~=~\int_{B_r(y^\nu)} \left(f\hbox{\bf 1}_\nu\right) (y)\left\{\int_{\R^n}e^{2\pi\i \left(\Phi(x,\xi)-y\cdot\xi\right)}\sigma(x,y,\xi)\phi^\nu(y)d\xi\right\}dy
\\\\ \ds
~=~\sum_{k=0}^\infty  \sum_{|\alpha|=k} \int_{\R^n}e^{2\pi\i\Phi(x,\xi)}\sigma_\alpha(x,\xi)\Hat{f}_\alpha(\xi)d\xi
\end{array}
\eeq
where $\sigma_\alpha(x,\xi)=\p_y^\alpha \sigma(x,y^\nu,\xi)\phi^\nu(y^\nu)/k!$ 
and $f_\alpha(y)= \left(f\hbox{\bf 1}_\nu\right)(y)\prod_{i=1}^n (y_i-y^\nu_i)^{\alpha_i}$  for $y\in B_{1/2n}(y^\nu)$ and $|\alpha|=\alpha_1+\alpha_2+\cdots+\alpha_n=k$. In particular, we have $\sum_{|\alpha|=k}\|f_\alpha\|_{\L^p(\R^n)}\leq\C 2^{-k}\|f\|_{\L^p(\R^n)}$. 

\begin{thm} Let $\sigma\in\S^0$ as (\ref{Class}). Fourier integral operator $\F$ defined in (\ref{Ff})-(\ref{nondegeneracy}) extends to a bounded operator 
\bel{L^2-result}
 \left\| \F f\right\|_{\L^2(\R^n)}~\leq~\C_{\sigma~\Phi}~\left\| f\right\|_{\L^2(\R^n)}.
 \eeq
\end{thm}

{\bf Proof:} First, it is suffice to consider $\mathcal{F}$ defined in (\ref{Ff'}). Furthermore, by Plancherel theorem, our estimates are reduced to a similar assertion for 
\bel{Sf}
\Big(\mathcal{S} f\Big)(x)~=~\int_{\R^n} e^{2\pi\i\Phi(x,\xi)} \sigma(x,\xi) f(\xi) d\xi
\eeq 
and its adjoint operator
\bel{S*f}
\Big(\mathcal{S}^* f\Big)(\xi)~=~\int_{\R^n} e^{-2\pi\i\Phi(x,\xi)} \bar{\sigma}(x,\xi) f(x) dx.
\eeq 
Let $\c$  be a small positive constant. We define an {\it narrow cone} as follows: whenever $\xi$ and $\eta$ belong to a same narrow cone and $|\eta|\leq|\xi|$, by writing $\eta=\rho\xi+\eta^\dagger$ for $0\leq\rho\leq1$ and $\eta^\dagger$ perpendicular to $\xi$, we require  $|\eta^\dagger|\leq\c\rho|\xi|$. The value of $\c$ depends on the phase function $\Phi$.
It is clear that  every $\mathcal{S}$ or $\mathcal{S}^*$ can be written as a finite sum of  partial operators,  whereas each one of them has  a symbol  supported on an narrow cone.

Recall the estimate given in {\it 3.1.1}, chapter IX of Stein \cite{Stein}.
We have
\bel{Phi_x est}
\left|\nabla_x \Big(\Phi(x,\xi)-\Phi(x,\eta)\Big)\right|~\ge~\C_\Phi~|\xi-\eta|
\eeq
whenever $\xi$ and $\eta$ belong to a same narrow cone.

From  direct computation, we have
\bel{S*Sf}
\begin{array}{cc}\ds
\Big(\mathcal{S}^*\mathcal{S} f\Big)(\xi)~=~\int_{\R^n} f(\eta)\mathfrak{S}^\sharp(\xi,\eta)d\eta,
\\\\ \ds
\mathfrak{S}^\sharp(\xi,\eta)~=~\int_{\R^n} e^{2\pi\i\left(\Phi(x,\eta)-\Phi(x,\xi)\right)} \sigma(x,\eta)\bar{\sigma}(x,\xi) dx.
\end{array}
\eeq
Since $\sigma(x,\xi)$ has a $x$-compact support, $\mathfrak{S}(\xi,\eta)$ in (\ref{S*Sf}) is uniformly bounded for $\sigma(x,\xi)\in\S^0$. On the other hand, we assume $\sigma(x,\xi)$ is supported on a narrow cone in the frequency space, with respect to a sufficiently small constant $\c$.

Recall from (\ref{Phi_x est}). 
An $N$-fold integration by parts  $w.r.t~x$ gives
\bel{Omega^sharp int by parts}
\begin{array}{lr}\ds
\left|\mathfrak{S}^\sharp(\xi,\eta)\right|~\leq~\C_{\Phi~N}~
\left|\xi-\eta\right|^{-N}
\left|\int_{\R^n} e^{2\pi\i\left(\Phi(x,\eta)-\Phi(x,\xi)\right)} \nabla_x^N\Big(\sigma(x,\eta)\bar{\sigma}(x,\xi)\Big)dx\right|
\end{array}
\eeq
for $\xi\neq\eta$ belong to a same narrow cone. 

The estimate in (\ref{Omega^sharp int by parts}) together with the differential inequality in  (\ref{Class}) imply
\bel{Omega^sharp est}
\left|\mathfrak{S}^\sharp(\xi,\eta)\right|~\leq~\C_{\sigma~\Phi~N}~\left({1\over 1+|\xi-\eta|}\right)^{N}
\eeq
for every $N\ge1$. 

By applying Minkowski integral inequality, we have
\bel{L^2 est S*S} 
\begin{array}{lr}\ds
\left\|\mathcal{S}^*\mathcal{S} f\right\|_{\L^2(\R^n)}~=~\left\{\int_{\R^n}\left|\int_{\R^n} f(\eta)\mathfrak{S}^\sharp(\xi,\eta)d\eta\right|^2d\xi\right\}^{1\over 2}
\\\\ \ds~~~~~~~~~~~~~~~~~~~~~~
~=~\left\{\int_{\R^n}\left|\int_{\R^n} f(\xi-\zeta)\mathfrak{S}^\sharp(\xi,\xi-\zeta)d\zeta\right|^2d\xi\right\}^{1\over 2}\qquad (~ \zeta=\xi-\eta ~) 
\\\\ \ds~~~~~~~~~~~~~~~~~~~~~~
~\leq~\C~\int_{\R^n} \left\{\int_{\R^n} \left|f(\xi-\zeta)\right|^2 \Big|\mathfrak{S}^\sharp(\xi,\xi-\zeta)\Big|^2 d\xi \right\}^{1\over 2 }d\zeta
\\\\ \ds~~~~~~~~~~~~~~~~~~~~~~
~\leq~\C_{\sigma~\Phi~N}~\int_{\R^n} \left\{\int_{\R^n} \left|f(\xi-\zeta)\right|^2 \left({1\over 1+|\zeta|}\right)^{2N} d\xi \right\}^{1\over 2 }d\zeta\qquad \hbox{\small{by (\ref{Omega^sharp est})}}
\\\\ \ds~~~~~~~~~~~~~~~~~~~~~~
~=~\C_{\sigma~\Phi~N}~\left\| f\right\|_{\L^2(\R^n)}\int_{\R^n}  \left({1\over 1+|\zeta|}\right)^{N}d\zeta
\\\\ \ds~~~~~~~~~~~~~~~~~~~~~~
~\leq~\C_{\sigma~\Phi}~\left\| f\right\|_{\L^2(\R^n)}
\end{array}
\eeq
provided that $N$ is sufficiently large.
\endproof
\begin{remark}
It can be easily seen that $\F_\ell$ defined in (\ref{delta_t})-(\ref{Partial}) satisfies (\ref{L^2-result}) as well. 
\end{remark}
Note that $\F_\ell$ can be written as a finite sum of
\bel{F_lf'}
\Big( \mathcal{F}_\ell f\Big)(x)~=~\int_{\R^n} e^{2\pi\i\Phi(x,\xi)}\sigma(x,\xi)\delta_\ell(\xi)\Hat{f}(\xi)d\xi
\eeq
whose $\L^2$-boundedness can be proved by carrying out same estimates in (\ref{Sf})-(\ref{L^2 est S*S}) with $\mathcal{S}f$ and $\mathcal{S}^*f$ replaced by
\bel{S_tf}
\Big(\mathcal{S}_\ell f\Big)(x)~=~\int_{\R^n} e^{2\pi\i\Phi(x,\xi)}\sigma(x,\xi)\delta_\ell(\xi)f(\xi)d\xi.
\eeq
and
\bel{S_t^*f}
\Big(\mathcal{S}^*_\ell f\Big)(\xi)~=~
\int_{\R^n} e^{-2\pi\i \Phi(x,\xi)} \bar{\sigma}(x,\xi)\bar{\delta}_\ell(\xi)f(x)dx
 \eeq
respectively.

\begin{thm} Let $\sigma\in\S^m$ for  $-n/2<m<0$. 
For every $\ell\ge0$, $\F_\ell$ defined in (\ref{delta_t})-(\ref{Partial}) extends to a bounded operator 
\bel{F_t 2,p result>}
\begin{array}{cc}\ds
 \left\| \F_\ell f\right\|_{\L^2(\R^n)}~\leq~\C_{p~\sigma~\Phi}~ 2^{\left({m\over n}\right) \ell}\left\| f\right\|_{\L^p(\R^n)}
\\\\ \ds
\hbox{for}\qquad{-m\over n}~=~{1\over p}-{1\over 2}
\end{array}
\eeq
and
\bel{F_t p',2 result>}
\begin{array}{cc}\ds
 \left\| \F_\ell f\right\|_{\L^{p\over p-1}(\R^n)}~\leq~\C_{p~\sigma~\Phi}~2^{\left({m\over n}\right) \ell}\left\| f\right\|_{\L^2(\R^n)}
 \\\\ \ds
 \hbox{for}\qquad
{-m\over n}~=~{1\over 2}-{p-1\over p}.
\end{array}
\eeq
\v
For every $\ell\leq0$, $\F_\ell$ defined in (\ref{delta_t})-(\ref{Partial}) extends to a bounded operator 
\bel{F_t 2,p result<}
\begin{array}{cc}\ds
 \left\| \F_\ell f\right\|_{\L^2(\R^n)}~\leq~\C_{p~\sigma~\Phi}~ 2^{-\left({m\over n}\right)(n-1) \ell}\left\| f\right\|_{\L^p(\R^n)}
\\\\ \ds
\hbox{for}\qquad{-m\over n}~=~{1\over p}-{1\over 2}
\end{array}
\eeq
and
\bel{F_t p',2 result<}
\begin{array}{cc}\ds
 \left\| \F_\ell f\right\|_{\L^{p\over p-1}(\R^n)}~\leq~\C_{p~\sigma~\Phi}~2^{-\left({m\over n}\right)(n-1) \ell}\left\| f\right\|_{\L^2(\R^n)}
 \\\\ \ds
 \hbox{for}\qquad
{-m\over n}~=~{1\over 2}-{p-1\over p}.
\end{array}
\eeq
\end{thm}

{\bf Proof:} From (\ref{F_lf'})-(\ref{S_t^*f}), it is suffice to estimate $\mathcal{F}_\ell$. 
Let $\ell\ge0$. We have
\bel{F_t decom}
\begin{array}{lr}\ds
\Big( \mathcal{F}_\ell f\Big)(x)~=~\int_{\R^n} e^{2\pi\i\Phi(x,\xi)}\sigma(x,\xi)\delta_\ell(\xi)\Hat{f}(\xi)d\xi
\\\\ \ds~~~~~~~~~~~~~~
~=~\int_{\R^n} e^{2\pi\i\Phi(x,\xi)}\sigma(x,\xi)|\xi|^{-m} \Big(\delta_\ell(\xi)\Hat{f}(\xi)|\xi|^m\Big)d\xi
\\\\ \ds~~~~~~~~~~~~~~
~=~ 2^{\left({m\over n}\right) \ell}\int_{\R^n} e^{2\pi\i\Phi(x,\xi)}\sigma(x,\xi)|\xi|^{-m} \left(\delta_\ell(\xi)\Hat{f}(\xi)2^{-\left({m\over n}\right) \ell}|\xi|^{m}\right)d\xi
\\\\ \ds~~~~~~~~~~~~~~
~\doteq~2^{\left({m\over n}\right) \ell}\int_{\R^n} e^{2\pi\i\Phi(x,\xi)}\sigma(x,\xi)|\xi|^{-m} \Big(\Hat{T_\ell f}\Big)(\xi)d\xi.
\end{array}
\eeq
By taking the inverse Fourier transform of $\Big(\Hat{T_\ell f}\Big)(\xi)$ defined implicitly in (\ref{F_t decom}), we have
\bel{convolution T_tf}
\begin{array}{cc}\ds
\Big(T_\ell f\Big)(x)~=~\int_{\R^n}f(y)\mathcal{K}_\ell(x-y)dy,
\\\\ \ds
\mathcal{K}_\ell(x)~=~\int_{\R^n} e^{2\pi\i x\cdot\xi} \delta_\ell(\xi) 2^{-\left({m\over n}\right)\ell}|\xi|^{m} d\xi.      
\end{array}     
\eeq
Let $\varphi$ be the smooth $bump$-function given in (\ref{varphi funda}). Define 
\bel{phi_j}
\phi_j(\xi)~=~\varphi(2^{-j}\xi)-\varphi(2^{-j+1}\xi),\qquad j\in\Z.
\eeq
From (\ref{convolution T_tf})-(\ref{phi_j}), we have
\bel{K_t sum}
\begin{array}{lr}\ds
\mathcal{K}_\ell(x)
~=~ 2^{-\left({m\over n}\right) \ell}\int_{\R^n} e^{2\pi\i x\cdot\xi} \delta_\ell(\xi) |\xi|^{m} d\xi
\\\\ \ds~~~~~~~~~
~=~\sum_{j\in\Z} 2^{-\left({m\over n}\right) \ell}\int_{\R^n} e^{2\pi\i x\cdot\xi} \delta_\ell(\xi) \phi_j(\xi)|\xi|^{m} d\xi
\\\\ \ds~~~~~~~~~
~=~\sum_{j\in\Z} 2^{-\left({m\over n}\right) \ell}\iint_{\R\times\R^{n-1}} e^{2\pi\i (z\tau+w\cdot\lambda)} \delta_\ell(\tau,\lambda) \phi_j(\tau,\lambda)\left(\tau^2+|\lambda|^2\right)^{m\over 2} d\tau d\lambda.
\end{array}
\eeq
Note that for $\xi$  in the support of $\phi_j(x)$, we have $2^{j-1}\leq|\xi|\leq2^{j+1},j\in\Z$. 

On the other hand, recall $\delta_\ell(\xi)$ defined in (\ref{delta_t}) supported on $\Lambda_\ell$ given in (\ref{Cone}).
We have $|\lambda|=\C 2^j$ and $|\tau|=\C 2^{j-\ell}$.  Observe that 
every $\partial_{\tau}$ acting on $\delta_\ell(\tau,\lambda) \phi_j(\tau,\lambda)\left(\tau^2+|\lambda|^2\right)^{m\over 2}$ gains a factor of $\C2^{-j+\ell}$ 
and every $\partial_{\lambda}$ acting on $\delta_\ell(\tau,\lambda) \phi_j(\tau,\lambda)\left(\tau^2+|\lambda|^2\right)^{m\over 2}$ gains a factor of $\C2^{-j}$.  

Moreover,  the volume of  $\supp\delta_\ell(\xi)\phi_j(\xi)$ is  bounded by $\C 2^{j-\ell}2^{j(n-1)}$.

An $N+M$-fold integration by parts $w.r.t~\xi=(\tau,\lambda)$ gives 
\bel{K_t by parts} 
\begin{array}{lr}\ds
\left| 2^{-\left({m\over n}\right) \ell}\int_{\R^n} e^{2\pi\i x\cdot\xi} \delta_\ell(\xi) \phi_j(\xi)|\xi|^{m} d\xi\right|
\\\\ \ds
~=~\left| 2^{-\left({m\over n}\right) \ell}\iint_{\R\times\R^{n-1}} e^{2\pi\i (z\tau+w\cdot\lambda)} \delta_\ell(\tau,\lambda) \phi_j(\tau,\lambda)\left(\tau^2+|\lambda|^2\right)^{m\over 2} d\tau d\lambda\right|
\\\\ \ds
~\leq~\C_{N~M}~ 2^{-\left({m\over n}\right) \ell}\left(2^{jm} 2^{j-\ell}2^{j(n-1)}\right)\left(2^{j-\ell}|z|\right)^{-N}\left(2^{j}|w|\right)^{-M}
\\\\ \ds
~=~\C_{N~M}~2^{(j-\ell)\left({n+m\over n}\right)}  \left(2^{j-\ell}|z|\right)^{-N} ~2^{j(n-1)\left({n+m\over n}\right)}\left(2^{j}|w|\right)^{-M}.
\end{array}
\eeq
We choose 
\bel{N,M tau}
\begin{array}{lr}\ds
N=0~~~\hbox{if}~~~|z|\leq2^{-j+\ell}\qquad \hbox{or}\qquad N=1 ~~~\hbox{if}~~~ |z|>2^{-j+\ell}; 
\\\\ \ds
M=0~~~ \hbox{if}~~~ |w|\leq2^{-j}\qquad \hbox{or}\qquad M=n-1~~~\hbox{if}~~~ |w|>2^{-j}.
\end{array}
\eeq
From (\ref{convolution T_tf}) and (\ref{K_t by parts}), we have
\bel{Sum K_t est} 
\begin{array}{lr}\ds
\left|\mathcal{K}_\ell(z,w)\right|~\leq~\sum_j \left| 2^{-\left({m\over n}\right) \ell}\iint_{\R\times\R^{n-1}} e^{2\pi\i (z\tau+w\cdot\lambda)} \delta_\ell(\tau,\lambda) \phi_j(\tau,\lambda)\left(\tau^2+|\lambda|^2\right)^{m\over 2} d\tau d\lambda\right|
\\\\ \ds~~~~~~~~~~~~~~~~
~\leq~\C_{N~ M} \sum_j 2^{(j-\ell)\left({n+m\over n}\right)}  \left(2^{j-\ell}|z|\right)^{-N} ~2^{j(n-1)\left({n+m\over n}\right)}\left(2^{j}|w|\right)^{-M}
\\\\ \ds~~~~~~~~~~~~~~~~ 
~\leq~\C_{N~M}  \left\{\sum_j2^{(j-\ell)\left({n+m\over n}\right)}  \left(2^{j-\ell}|z|\right)^{-N}\right\}
\left\{\sum_j 2^{j(n-1)\left({n+m\over n}\right)}\left(2^{j}|w|\right)^{-M}\right\}
\\\\ \ds~~~~~~~~~~~~~~~~
~=~\C~  \left\{\sum_{|z|\leq2^{-j+\ell}} 2^{(j-\ell)\left({n+m\over n}\right)}~+~\sum_{|z|>2^{-j+\ell}}2^{(j-\ell)\left({n+m\over n}\right)} \left(2^{j-\ell}|z|\right)^{-1}\right\}
\\\\ \ds~~~~~~~~~~~~~~~~~~~~~~~~~
\left\{\sum_{|w|\leq2^{-j}} 2^{j(n-1)\left({n+m\over n}\right)}~+~\sum_{|w|>2^{-j}}2^{j(n-1)\left({n+m\over n}\right)} \left(2^{j}|w|\right)^{-(n-1)}\right\}\qquad\hbox{\small{by (\ref{N,M tau})}}
\\\\ \ds~~~~~~~~~~~~~~~~
~\leq~\C~ \left\{\left({1\over |z|}\right)^{n+m\over n}+\left({1\over |z|}\right) \sum_{|z|>2^{-j+\ell}} 2^{(j-\ell)\left({m\over n}\right)}\right\}\left\{\left({1\over |w|}\right)^{(n-1)\left({n+m\over n}\right)}+\left({1\over |w|}\right)^{n-1} \sum_{|w|>2^{-j}} 2^{j(n-1)\left({m\over n}\right)}\right\}
\\ \ds~~~~~~~~~~~~~~~~~~~~~~~~~~~~~~~~~~~~~~~~~~~~~~~~~~~~~~~~~~~~~~~~~~~~~~~~~~~~~~~~~~~~~~~~~~~~~~~~~~~~~~~~~~~~~~~~~~~~~~~~~~~~~~~~~~~~~~~~~
 \hbox{\small{( $m<0$ )}}
\\ \ds~~~~~~~~~~~~~~~~
~\leq~\C~ \left({1\over |z|}\right)^{n+m\over n}\left({1\over|w|}\right)^{(n-1)\left({n+m\over n}\right)}.
\end{array}
\eeq

Let $\alpha={-m\over n},\beta=-m\left({n-1\over n}\right)$ of which ${\alpha}={\beta\over n-1}={1\over p}-{1\over 2}$. By applying Hardy-Littlewood-Sobolev theorem \cite{Hardy-Littlewood}-\cite{Sobolev} on the subspaces $\R$ and $\R^{n-1}$, we have
\bel{T_t regularity}
\begin{array}{lr}\ds
\left\|T_\ell f\right\|_{\L^2(\R^n)}~=~\left\{\iint_{\R\times\R^{n-1}} \left|\iint_{\R\times\R^{n-1}} f(u,v)\mathcal{K}_\ell(z-u,w-v)dudv\right|^2dzdw\right\}^{1\over 2}
\\\\ \ds
~\leq~\C~\left\{\iint_{\R\times\R^{n-1}} \left\{\iint_{\R\times\R^{n-1}} \left|f(u,v)\right|\left({1\over|z-u|}\right)^{n+m\over n}\left({1\over|w-v|}\right)^{(n-1)\left({n+m\over n}\right)}dudv\right\}^2dzdw\right\}^{1\over 2}~~ \hbox{\small{by (\ref{Sum K_t est})}}
\\\\ \ds
~\leq~\C~\left\{\int_{\R^{n-1}} \left\{\int_{\R}\left\{\int_{\R^{n-1}} \left|f(z,v)\right|\left({1\over|w-v|}\right)^{(n-1)\left({n+m\over n}\right)}dv\right\}^pdz\right\}^{2\over p}dw\right\}^{1\over 2}
\\\\ \ds
~\leq~\C~\left\{\int_{\R^{n}} \left\{\int_{\R^{N-1}}\left\{\int_{\R^{n-1}} \left|f(z,v)\right|\left({1\over|w-v|}\right)^{(n-1)\left({n+m\over n}\right)}dv\right\}^2dw\right\}^{p\over 2}dz\right\}^{1\over p}
\\ \ds~~~~~~~~~~~~~~~~~~~~~~~~~~~~~~~~~~~~~~~~~~~~~~~~~~~~~~~~~~~~~~~~~~~~~~~~~~~~~
\hbox{\small{by Minkowski integral inequality}}
\\ \ds
~\leq~\C_p~\left\| f\right\|_{\L^p(\R^n)}.
\end{array}
\eeq
Note that  $\sigma(x,\xi)|\xi|^{-m}\in\S^0$ in (\ref{F_t decom}). By applying {\bf Theorem 2.1} and  using (\ref{T_t regularity}), we have
\bel{F_t f regularity}
\begin{array}{lr}\ds
\left\|\mathcal{F}_\ell f\right\|_{\L^2(\R^n)}~\leq~\C_{\sigma~\Phi}~ 2^{\left({m\over n}\right)\ell} \left\| T_\ell f\right\|_{\L^2(\R^n)}
\\\\ \ds~~~~~~~~~~~~~~~~~~~
~\leq~\C_{p~\sigma~\Phi}~2^{\left({m\over n}\right)\ell} \left\|f\right\|_{\L^p(\R^n)}\qquad
\hbox{for}\qquad
{-m\over n}~=~{1\over p}-{1\over 2}.
\end{array}
\eeq
From direct computation, the adjoint operator
\bel{F^*_t f}
\begin{array}{lr}\ds
\Big(\mathcal{F}^*_\ell f\Big)(x)~=~\int_{\R^n} f(y) \left\{\int_{\R^n} e^{2\pi\i \left(x\cdot\xi-\Phi(y,\xi)\right)}\bar{\sigma}(y,\xi)\bar{\delta}_\ell(\xi)d\xi\right\}dy
\\\\ \ds~~~~~~~~~~~~~~
~=~\int_{\R^n} e^{2\pi\i x\cdot \xi}  \left\{ \int_{\R^n} e^{-2\pi\i \Phi(y,\xi)} \bar{\sigma}(y,\xi)\bar{\delta}_\ell(\xi)f(y)dy\right\}d\xi
\\\\ \ds~~~~~~~~~~~~~~
~=~\int_{\R^n} e^{2\pi\i x\cdot\xi} \Big(\mathcal{S}^*_\ell f\Big)(\xi)d\xi.
\end{array}
\eeq
We aim to show $\F_\ell\colon\L^2(\R^n)\mt\L^{p\over p-1}(\R^n)$ for ${-m\over n}={1\over 2}-{p-1\over p}$ by proving  $\F^*_\ell\colon\L^p(\R^n)\mt\L^2(\R^n)$ for ${-m\over n}={1\over p}-{1\over 2}$, with the desired operator norm for every $\ell$. By Plancherel theorem, we reduce our assertion for $\mathcal{S}^*_\ell$ given in (\ref{S_t^*f}).    
Observe that 
\bel{S^*_t L^2-norm}
\begin{array}{lr}\ds
\left\| \mathcal{S}^*_\ell f\right\|_{\L^2(\R^n)}^2~=~\int_{\R^n} \Big(\mathcal{S}_\ell\mathcal{S}^*_\ell f\Big)(x) f(x)dx
\\\\ \ds~~~~~~~~~~~~~~~~~~~
~\leq~\left\|\mathcal{S}_\ell\mathcal{S}^*_\ell f\right\|_{\L^{p\over p-1}(\R^n)}\left\| f\right\|_{\L^p(\R^n)}\qquad \hbox{\small{by H\"{o}lder inequality.}}
\end{array}
\eeq
We prove $\mathcal{S}^*_\ell\colon\L^p(\R^n)\mt\L^2(\R^n)$ for ${-m\over n}={1\over p}-{1\over 2}$ by showing $\mathcal{S}_\ell\mathcal{S}^*_\ell\colon\L^p(\R^n)\mt\L^{p\over p-1}(\R^n)$ for $-{2m\over n}={1\over p}-{p-1\over p}$, with the desired operator norm for every $\ell\ge0$.

From direct computation, we have
\bel{S_tS^*_t}
\begin{array}{cc}\ds
\Big(\mathcal{S}_\ell\mathcal{S}^*_\ell f\Big)(x)~=~\int_{\R^n} f(y)\mathfrak{S}^\flat_\ell(x,y)dy,
\\\\ \ds
\mathfrak{S}^\flat_\ell(x,y)~=~\int_{\R^n} e^{2\pi\i\left(\Phi(x,\xi)-\Phi(y,\xi)\right)} \sigma(x,\xi)\delta_\ell(\xi)\bar{\sigma}(y,\xi)\bar{\delta}_\ell(\xi)d\xi.
\end{array}
\eeq
Note that
\bel{nabla Phi}
\nabla_\xi \left(\Phi(x,\xi)-\Phi(y,\xi)\right)~=~\Phi_{x\xi}(x,\xi)(x-y)+\O\left(|x-y|^2\right).
\eeq
We momentarily assume that $\sigma\in\S^m$ has a sufficiently small  support in $x$.

From (\ref{nabla Phi}) and keeping in mind that $\Phi(x,\xi)$ satisfies the nondegeneracy condition in (\ref{nondegeneracy}), we have
\bel{nabla Phi size}
\left|\nabla_\xi \left(\Phi(x,\xi)-\Phi(y,\xi)\right)\right|~\ge~\C_{\Phi}~|x-y|.
\eeq
By changing dilations
\bel{t dila}
\begin{array}{cc}\ds
x~=~\mathfrak{L}^{-1}x'~\doteq~\left(2^{\ell}z', w\right),\qquad y~=~\mathfrak{L}^{-1}y'~\doteq~\left(2^{\ell}u',v\right)
\\\\ \ds
\hbox{and}\qquad \xi~=~\mathfrak{L}\xi'~\doteq~\left(2^{-\ell}\tau', \lambda\right).
\end{array}
\eeq
Observe that $\delta_\ell(\mathfrak{L}\xi')=\delta_o(\xi')$ by definition of $\delta_\ell(\xi)$ in (\ref{delta_t}).
Let $\phi_j(\xi)$ be defined in (\ref{phi_j}). We write
\bel{Sum S^flat dila}
\begin{array}{lr}\ds
\mathfrak{S}^\flat_\ell\left(x,y\right)~=~\mathfrak{S}^\flat_\ell\left(\mathfrak{L}^{-1}x',\mathfrak{L}^{-1}y'\right)
\\\\ \ds
~=~2^{-\ell}\int_{\R^n} e^{2\pi\i\left(\Phi(\mathfrak{L}^{-1}x',\mathfrak{L}\xi')-\Phi(\mathfrak{L}^{-1}y',\mathfrak{L}\xi')\right)} \sigma(\mathfrak{L}^{-1}x',\mathfrak{L}\xi')\delta_o(\xi')\bar{\sigma}(\mathfrak{L}^{-1}y',\mathfrak{L}\xi')\bar{\delta}_o(\xi')d\xi'
\\\\ \ds
~=~\sum_{j\in\Z} 2^{-\ell}\int_{\R^n} e^{2\pi\i\left(\Phi(\mathfrak{L}^{-1}x',\mathfrak{L}\xi')-\Phi(\mathfrak{L}^{-1}y',\mathfrak{L}\xi')\right)} \sigma(\mathfrak{L}^{-1}x',\mathfrak{L}\xi')\delta_o(\xi')\bar{\sigma}(\mathfrak{L}^{-1}y',\mathfrak{L}\xi')\bar{\delta}_o(\xi')\phi_j(\xi')d\xi'.
\end{array}
\eeq
Recall  $2^{j-1}\leq|\xi'|\leq2^{j+1},j\in\Z$ for $\xi'$  in the support of $\phi_j(\xi')$. Observe that   every $\partial_{\xi'}$ acting on $ \sigma(\mathfrak{L}^{-1}x',\mathfrak{L}\xi')\delta_o(\xi')\bar{\sigma}(\mathfrak{L}^{-1}y',\mathfrak{L}\xi')\bar{\delta}_o(\xi')\phi_j(\xi')$ gains a factor of $\C2^{-j}$. 

Let $\Phi_{x\xi}(x,\xi)$ denote the $n\times n$-matrix $\left({\p^2\Phi\over\p x_i\p\xi_j}\right)(x,\xi)$ in (\ref{nondegeneracy}). We have
\bel{det Phi_xxi}
\det\Phi_{x\xi}(x,\xi)~=~\det\Phi_{x'\xi'}(\mathfrak{L}^{-1}x',\mathfrak{L}\xi').
\eeq
Indeed,    $2^{\ell}$ appears at  the first column of $\Phi_{x'\xi'}(\mathfrak{L}^{-1}x',\mathfrak{L}\xi')$ and  $2^{-\ell}$ appears at  the first row of $\Phi_{x'\xi'}(\mathfrak{L}^{-1}x',\mathfrak{L}\xi')$ respectively. 

From (\ref{nabla Phi size}), we thus have
\bel{nabla Phi size new}
\left|\nabla_{\xi'} \left(\Phi(\mathfrak{L}^{-1}x',\mathfrak{L}\xi')-\Phi(\mathfrak{L}^{-1}y',\mathfrak{L}\xi')\right)\right|~\ge~\C_{\Phi}~|x'-y'|.
\eeq

Note that $\sigma\in\S^m$ and
the support of $\phi_j(\xi)$ has a volume bounded by $\C2^{jn}$.
An $N+M$-fold integration by parts $w.r.t~ \xi$ gives
\bel{S^flat by parts}
\begin{array}{lr}\ds
 2^{-\ell}\left|\int_{\R^n} e^{2\pi\i\left(\Phi(\mathfrak{L}^{-1}x,\mathfrak{L}\xi)-\Phi(\mathfrak{L}^{-1}y,\mathfrak{L}\xi)\right)} \sigma(\mathfrak{L}^{-1}x,\mathfrak{L}\xi)\delta_o(\xi)\bar{\sigma}(\mathfrak{L}^{-1}y,\mathfrak{L}\xi)\bar{\delta}_o(\xi)\phi_j(\xi)d\xi\right|
\\\\ \ds
~\leq~\C_\Phi~ 2^{-\ell}\left(2^{2mj} 2^{jn}\right)\left(2^j |x'-y'|\right)^{-N-M}
\\\\ \ds
~\leq~\C_\Phi~\left( 2^{\left({2m\over n}\right)j} 2^{(n-1)\left({2m\over n}\right)j}\right)\left( 2^{j-\ell}2^{j(n-1)}\right) \left( 2^j |z'-u'|\right)^{-N}\left( 2^j |w-v|\right)^{-M}
\\\\ \ds
~\leq~\C_\Phi~ 2^{\left({2m\over n}\right)\ell}\left[2^{(j-\ell)\left({n+2m\over n}\right)}  \left( 2^j |z'-u'|\right)^{-N}\right]\left[2^{j(n-1)\left({n+2m\over n}\right)}  \left( 2^j |w-v|\right)^{-M}\right]
\\\\ \ds
~\leq~\C_\Phi~ 2^{\left({2m\over n}\right)\ell}\left[2^{(j-\ell)\left({n+2m\over n}\right)}  \left( 2^{j-\ell} |z-u|\right)^{-N}\right]\left[2^{j(n-1)\left({n+2m\over n}\right)}  \left( 2^j |w-v|\right)^{-M}\right].
\end{array}
\eeq
We choose 
\bel{N,M lambda}
\begin{array}{lr}\ds
N=0~~~\hbox{if}~~~ |z-u|\leq2^{-j+\ell}\qquad \hbox{or}\qquad N=1~~~\hbox{if}~~~ |z-u|>2^{-j+\ell},
\\\\ \ds
M=0~~~\hbox{if}~~~ |w-v|\leq2^{-j}\qquad \hbox{or}\qquad M=n-1~~~\hbox{if}~~~ |w-v|>2^{-j}.
\end{array}
\eeq
From (\ref{Sum S^flat dila}) and (\ref{S^flat by parts}), we have
\bel{Sum S^flat est} 
\begin{array}{lr}\ds
\left|\mathfrak{S}^\flat_\ell\left(x,y\right)\right|~\leq~\C_{\Phi~N~M}~ \sum_{j\in\Z} 2^{\left({2m\over n}\right)\ell}\left[2^{(j-\ell)\left({n+2m\over n}\right)}  \left( 2^{j-\ell} |z-u|\right)^{-N}\right]\left[ 2^{j(n-1)\left({n+2m\over n}\right)}  \left( 2^j |w-v|\right)^{-M}\right]
\\\\ \ds~~~~~~~~~~~~~~~ 
~\leq~\C_{\Phi~N~M}~   2^{\left({2m\over n}\right)\ell} \left\{ \sum_{j\in\Z} 2^{(j-\ell)\left({n+2m\over n}\right)}  \left( 2^{j-\ell} |z-u|\right)^{-N} \right\}\left\{\sum_{j\in\Z} 2^{j(n-1)\left({n+2m\over n}\right)}  \left( 2^j |w-v|\right)^{-M}\right\}
\\\\ \ds~~~~~~~~~~~~~~~
~=~\C_{\Phi}~ 2^{\left({2m\over n}\right)\ell} \left\{\sum_{|z-u|\leq2^{-j+\ell}} 2^{(j-\ell)\left({n+2m\over n}\right)}~+~\sum_{|z-u|>2^{-j+\ell}}2^{(j-\ell)\left({n+2m\over n}\right)} \left(2^{j-\ell}|z-u|\right)^{-1}\right\}
\\\\ \ds~~~~~~~~~~~~~~~~~~~~~~~~~~~~~~~~~~~~~~
\left\{\sum_{|w-v|\leq2^{-j}} 2^{j(n-1)\left({n+2m\over n}\right)}~+~\sum_{|w-v|>2^{-j}} 2^{j(n-1)\left({n+2m\over n}\right)} \left(2^{j}|w-v|\right)^{-(n-1)}\right\}\qquad\hbox{\small{by (\ref{N,M lambda})}}
\\\\ \ds~~~~~~~~~~~~~~~
~\leq~\C_\Phi~ 2^{\left({2m\over n}\right)\ell} \left\{ \left({1\over |z-u|}\right)^{n+2m\over n}~+~\left({1\over |z-u|}\right) \sum_{|z-u|>2^{-j+\ell}} 2^{(j-\ell)\left({2m\over n}\right)}\right\}
\\\\ \ds~~~~~~~~~~~~~~~~~~~~~~~~~~~~~~~~~~~~~~
 \left\{ \left({1\over |w-v|}\right)^{(n-1)\left({n+2m\over n}\right)}~+~\left({1\over |w-v|}\right)^{n-1} \sum_{|w-v|>2^{-j}} 2^{j(n-1)\left({2m\over n}\right)}\right\}
\qquad \hbox{\small{($m<0$)}}
\\\\ \ds~~~~~~~~~~~~~~~
~\leq~\C_\Phi~ 2^{\left({2m\over n}\right)\ell}\left({1\over |z-u|}\right)^{n+2m\over n}\left({1\over |w-v|}\right)^{(n-1)\left({n+2m\over n}\right)}.
\end{array}
\eeq
Recall that $\sigma(x,\xi)$ is assumed to have a sufficiently small $x$-support,  through (\ref{nabla Phi size})-(\ref{Sum S^flat est}). Because $\sigma(x,\xi)$ has a compact support in $x$,  it can be written as a finite sum of symbols having this extra restriction.  

Let ${-2m\over n}={1\over p}-{p-1\over p}$. Recall from (\ref{S_tS^*_t}) and the estimate in (\ref{Sum S^flat est}). By carrying out the iteration argument given in (\ref{T_t regularity}) and taking into account that the implied constants depend also on the size of the $x$-support of $\sigma(x,\xi)$,
we have 
\bel{S_tS^*_t regularity}
\begin{array}{lr}\ds
\left\|\mathcal{S}_\ell\mathcal{S}^*_\ell f\right\|_{\L^{p\over p-1}(\R^n)}~=~\left\{\int_{\R^n}\left|\int_{\R^n} f(y)\mathfrak{S}^\flat_\ell(x,y)dy\right|^{p\over p-1}dx\right\}^{p-1\over p}
\\\\ \ds
~\leq~\C_{\sigma~\Phi}~ 2^{\left({2m\over n}\right)\ell}\left\{\iint_{\R\times\R^{n-1}}\left\{\iint_{\R\times\R^{n-1}} |f(u,v)|\left({1\over |z-u|}\right)^{n+2m\over n}\left({1\over |w-v|}\right)^{(n-1)\left({n+2m\over n}\right)}dudv\right\}^{p\over p-1}dzdw\right\}^{p-1\over p}
\\\\ \ds
~\leq~\C_{p~\sigma~\Phi}~ 2^{\left({2m\over n}\right)\ell}~\left\| f\right\|_{\L^p(\R^n)}.
\end{array}
\eeq
By putting together (\ref{S^*_t L^2-norm}) and (\ref{S_tS^*_t regularity}), we find
$\left\|\mathcal{S}^*_\ell f\right\|_{\L^{2}(\R^n)}\leq\C_{p~\sigma~\Phi}~ 2^{\left({m\over n}\right)\ell}~\left\| f\right\|_{\L^p(\R^n)}$
as  desired.

For $\ell\leq0$, it is equivalent to consider 
\bel{delta_l}
\deltaup_\ell(\xi)~\doteq~\varphi\left(2^{\ell}{\lambda\over|\tau|}\right)-\varphi\left(2^{\ell+1}{\lambda\over|\tau|}\right),\qquad \ell\ge0.
\eeq
The result in (\ref{F_t 2,p result<})-(\ref{F_t p',2 result<}) can be proved in the same estimates as (\ref{F_t decom})-(\ref{S_tS^*_t regularity})  by switching the role of $\tau$ and $\lambda$ together with their dual variables $(z,w)$ and $(u,v)$ respectively.\endproof

\section{An heuristic argument}
\setcounter{equation}{0}
 Let $B_r(x_o)\subset\R^{n}$ be a ball centered on 
$x_o\in\R^n$ with radius $r>0$ and 
 $a$ denote an $\H^1$-{\it atom} associated to $B_r(x_o)$.
We aim to show that for $\sigma\in\S^{-{n-1\over 2}}$, 
\bel{H^1 est F_ta}
\begin{array}{lr}\ds
\int_{\R^n} \left|\Big(\F_\ell a\Big)(x)\right|dx~\leq~\C_{\sigma~\Phi}\left\{\begin{array}{lr}\ds 2^{-\left({n-1\over 2n}\right)\ell},\qquad \ell\ge0,
\\\\ \ds
2^{\left(n-1\right)\left({n-1\over 2n}\right)\ell},\qquad \ell\leq0.
\end{array}\right.
\end{array}
\eeq 
In the following, we focus on $\ell\ge0$ whereas the estimates for $\ell\leq0$ are handled similarly.  

We consider the {\it region of influence}, denoted by $\Q_r=\Q_r(x_o)$, satisfying
\bel{Q norm}
|\Q_r(x_o)|~\leq~\C_\sigma~r.
\eeq
The actual set of $\Q_r(x_o)$ 
will be explicitly constructed in section 6.

By applying Schwartz inequality, we have
\bel{Local est}
\int_{\Q_r(x_o)} \left|\Big(\F_\ell a\Big)(x)\right|dx~\leq~|\Q_r(x_o)|^{1\over 2}\left\|\F_\ell a\right\|_{\L^2(\R^n)}~\leq~\C_\sigma~r^{1\over 2}\left\|\F_\ell a\right\|_{\L^2(\R^n)}.
\eeq
On the other hand, the $\H^1$-atom $a$ associated to the ball $B_r(x_o)$ satisfies  
\bel{H^1 atom size}
\|a\|_{\L^p(\R^n)}~=~\left\{\int_{B_r(x_o)} \left|a(x)\right|^p dx \right\}^{1\over p} ~\leq~ |B_r(x_o)|^{-1+{1\over p}},\qquad 1\leq p<\infty.
\eeq
Let $\sigma\in\S^{-{n-1\over 2}}$. By applying {\bf Theorem 3.1}, the estimate in (\ref{F_t 2,p result>}) implies
\bel{F_t a Result p 2}
\begin{array}{cc}\ds
\left\|\F_\ell a\right\|_{\L^2(\R^n)}~\leq~\C_{p~\sigma~\Phi}~2^{-\left({n-1\over 2n}\right)\ell}~ \|a\|_{\L^p(\R^n)}
\\\\ \ds
 \hbox{for}\qquad 
{1\over p}~=~{1\over 2}~+~{n-1\over 2n}.
\end{array}
\eeq
From (\ref{Local est})-(\ref{F_t a Result p 2}), we have
\bel{Local est 2}
\begin{array}{lr}\ds
\int_{\Q_r(x_o)} \left|\Big(\F_\ell a\Big)(x)\right|dx~\leq~\C_{p~\sigma~\Phi}~2^{-\left({n-1\over 2n}\right)\ell}~ \|a\|_{\L^p(\R^n)}
\\\\ \ds~~~~~~~~~~~~~~~~~~~~~~~~~~~~~~~
~\leq~\C_{p~\sigma~\Phi}~ r^{1\over 2}r^{n\left(-1+{1\over p}\right)}~2^{-\left({n-1\over 2n}\right)\ell}
\\\\ \ds~~~~~~~~~~~~~~~~~~~~~~~~~~~~~~~
~=~\C_{\sigma~\Phi}~2^{-\left({n-1\over 2n}\right)\ell}
\end{array}
\eeq
where $-1+{1\over p}={-1\over 2}+{n-1\over 2n}={-1\over 2n}$.
Suppose that we can show
\bel{Comple est}
\int_{{^c}\Q_r(x_o)} \left|\Big(\F_\ell a\Big)(x)\right|dx~\leq~\C_{\sigma~\Phi}~2^{-\left({n-1\over 2n}\right)\ell}.
\eeq
Together with (\ref{Local est 2}), we obtain (\ref{H^1 est F_ta}) for $\ell\ge0$.
Recall the characterization of Hardy spaces  given  by Fefferman and Stein \cite{FC.S}. For $f\in\H^1\left(\R^n\right)$, it  can be written as $\sum_{k=1}^\infty c_k a_k$ 
of which every  $a_k$ is an $\H^1$-atom  and $\sum_{k=1}^\infty |c_k|\leq\C$. 
Therefore,  (\ref{H^1 est F_ta}) further  implies
\bel{H^1 est}
\int_{\R^n} \left|\Big(\F_\ell f\Big)(x)\right|dx~\leq~\C_{\sigma~\Phi}~ 2^{-\left({n-1\over 2n}\right)\ell}\left\| f\right\|_{\H^1(\R^n)}.
\eeq
Recall from (\ref{Partial}). Consider the adjoint operator  
\bel{F*_t singular}
\begin{array}{cc}\ds
\Big(\F^*_\ell f\Big)(x)~=~\int_{\R^n} f(y)\Omega^*_\ell(x,y)dy,
\\\\ \ds
 \Omega^*_\ell(x,y)~=~\int_{\R^n} e^{2\pi\i\left(x\cdot\xi-\Phi(y,\xi)\right)}\bar{\delta}_\ell(\xi)\bar{\sigma}(x,y,\xi)d\xi.
 \end{array}
 \eeq 
By duality between $\H^1$ and $\BMO$ spaces, as investigated by Fefferman \cite{Fefferman}, (\ref{H^1 est}) implies
\bel{L^infty Est}
\left\|\F^*_\ell f\right\|_{\B\M\O(\R^n)}~\leq~\C_{\sigma~\Phi}~2^{-\left({n-1\over 2n}\right)\ell}~\left\|f\right\|_{\L^\infty(\R^n)}.
\eeq
The {\it region of influence} associated to  $\F^*_\ell$,  denoted by $\Q^*_r(x_o)$,  satisfies 
\bel{Q* norm}
\left| \Q^*_r(x_o)\right|~\leq~\C~r.
\eeq
The actual set of $\Q^*_r(x_o)$ will  be defined explicitly in section 6.

By  using the estimate in (\ref{F_t p',2 result>}), we have
\bel{F_t a Result 2 p' and F^*_t a Result p 2}
\begin{array}{ccc}\ds
\left\|\F_\ell a\right\|_{\L^2(\R^n)}~\leq~\C_{p~\sigma~\Phi}~2^{-\left({n-1\over 2n}\right)\ell}~ \|a\|_{\L^p(\R^n)}
\\\\ \ds
 \hbox{for}\qquad 
{1\over 2}~=~{p-1\over p}~+~{n-1\over 2n}
\\\\ \ds
\Longleftrightarrow\qquad\left\|\F^*_\ell a\right\|_{\L^2(\R^n)}~\leq~\C_{p~\sigma~\Phi}~2^{-\left({n-1\over 2n}\right)\ell}~ \|a\|_{\L^p(\R^n)}
\\\\ \ds
\hbox{for}\qquad 
{1\over p}~=~{1\over 2}~+~{n-1\over 2n}.
\end{array}
\eeq
By carrying out same estimates  in (\ref{Local est})-(\ref{Comple est}), with $\F_\ell$ and $\Q_r(x_o)$ replaced by $\F^*_\ell$ and $\Q^*_r(x_o)$ respectively, we have
\bel{H^1 est *}
\int_{\R^n} \left|\Big(\F^*_\ell f\Big)(x)\right|dx~\leq~\C_{\sigma~\Phi}~2^{-\left({n-1\over 2n}\right)\ell}~\left\| f\right\|_{\H^1(\R^n)}.
\eeq
Hence that the duality between $\H^1$ and $\BMO$ spaces implies
\bel{L^infty Est*}
\left\|\F_\ell f\right\|_{\B\M\O(\R^n)}~\leq~\C_{\sigma~\Phi}~2^{-\left({n-1\over 2n}\right)\ell}~\left\|f\right\|_{\L^\infty(\R^n)}.
\eeq
\v

We now proceed to an interpolation argument set out in {\bf 5.2} chapter IV of \cite{Stein}. 

Consider an analytic family of operators $\F_{\ell~\hbox{z}}$ defined on the strip $\left\{\hbox{z}\in\mathbb{C}\colon0<\Re(\hbox{z})<1\right\}$ by 
\bel{F_z}
\Big(\F_{\ell~\hbox{z}} f\Big)(x)~=~e^{(\hbox{z}-\vartheta)^2}\iint_{\R^n\times\R^n} e^{2\pi\i \left(\Phi(x,\xi)-y\cdot\xi\right)} \sigma(x,y,\xi)\left(1+|\xi|^2\right)^{\gamma(\hbox{\small{z}})\over 2}\delta_\ell(\xi) f(y)dyd\xi
\eeq
where $\sigma\in\S^m$ and
\bel{gamma_z, theta}
\gamma(\hbox{z})~=~-m-{\hbox{z}(n-1)\over 2},\qquad \vartheta~=~-{2m\over n-1}.
\eeq
Note that $e^{(\hbox{z}-\vartheta)^2}$ decays rapidly as $|\Im(\hbox{z})|\mt\infty$.

For every $\hbox{z}$ in the strip, we have $e^{(\hbox{z}-\vartheta)^2}\sigma(x,y,\xi)\left(1+|\xi|^2\right)^{\gamma(\hbox{\small{z}})\over 2}\in\S^0$.  {\bf Remark 3.1} implies
\bel{L^2 Est z}
\left\|\F_{\ell~\hbox{z}} f\right\|_{\L^2(\R^n)} ~\leq~\C_{\sigma~\Phi} \left\|f\right\|_{\L^2(\R^n)},~~0~<~\Re(\hbox{z})~<~1,~~ -\infty~<~\Im(\hbox{z})~<~\infty.
\eeq 
When $\Re(\hbox{z})=0$, we have $e^{(\hbox{z}-\vartheta)^2}\sigma(x,y,\xi)\left(1+|\xi|^2\right)^{\gamma(\hbox{\small{z}})\over 2}\in\S^0$ and
\bel{L^2 Est t}
\left\|\F_{\ell~\Im(\hbox{z})\i} f\right\|_{\L^2(\R^n)} ~\leq~\C_{\sigma~\Phi} \left\|f\right\|_{\L^2(\R^n)},\qquad -\infty~<~\Im(\hbox{z})~<~\infty.
\eeq 
When $\Re(\hbox{z})=1$, we have $e^{(\hbox{z}-\vartheta)^2}\sigma(x,y,\xi)\left(1+|\xi|^2\right)^{\gamma(\hbox{\small{z}})\over 2}\in\S^{-{n-1\over 2}}$.  The estimate in (\ref{L^infty Est*}) implies
\bel{L^infty Est t}
\left\|\F_{\ell~1+ \Im(\hbox{z})\i} f\right\|_{\B\M\O(\R^n)}~\leq~\C_{\sigma~\Phi}~2^{-\left({n-1\over 2n}\right)\ell}~\left\|f\right\|_{\L^\infty(\R^n)},\qquad -\infty~<~\Im(\hbox{z})~<~\infty.
\eeq
By applying the complex interpolation theorem of Fefferman and Stein \cite{FC.S}, we obtain 
\bel{L^p Result vartheta}
\begin{array}{lr}\ds
\left\|\F_{\ell~\vartheta} f\right\|_{\L^p(\R^n)}~\leq~ \C_{p~\sigma~\Phi}~2^{-\vartheta\left({n-1\over 2n}\right)\ell}~\left\|f\right\|_{\L^p(\R^n)},\qquad {1\over p}~=~{1\over 2}(1-\vartheta)
\\\\ \ds~~~~~~~~~~~~~~~~~~~~~~
~=~ \C_{p~\sigma~\Phi}~2^{\left({m\over n}\right)\ell}~\left\|f\right\|_{\L^p(\R^n)}
\end{array}
\eeq
where $\vartheta=-{2m\over n-1}$ from (\ref{gamma_z, theta}).
Observe that $\F_{\ell~\vartheta}=\F_\ell$ and ${1\over2}-{1\over p}={-m\over n-1}$.

On the other hand, consider  $\F^*_{\ell~\hbox{z}}$ defined on the strip $\left\{\hbox{z}\in\mathbb{C}\colon0<\Re(\hbox{z})<1\right\}$ by 
\bel{F^*_z}
\Big(\F^*_{\ell~\hbox{z}} f\Big)(x)~=~e^{(\hbox{z}-\vartheta)^2}\iint_{\R^n\times\R^n} e^{2\pi\i \left(x\cdot\xi-\Phi(y,\xi)\right)} \bar{\sigma}(x,y,\xi)\left(1+|\xi|^2\right)^{\gamma(\hbox{\small{z}})\over 2}\bar{\delta}_\ell(\xi) f(y)dyd\xi
\eeq
where $\sigma\in\S^m$ and $\gamma(\hbox{z}), \vartheta$ are defined in (\ref{gamma_z, theta}).

Note that same estimates hold in (\ref{L^2 Est z})-(\ref{L^infty Est t}) for $\F^*_{\ell~\hbox{z}}$. In particular,  we use  (\ref{L^infty Est}) instead of (\ref{L^infty Est*}) to show that $\F^*_{\ell~\hbox{z}}$ satisfies the norm inequality in (\ref{L^infty Est t}). 

By applying the desired complex interpolation, we obtain 
\bel{L^p Result vartheta *}
\left\|\F^*_{\ell~\vartheta} f\right\|_{\L^p(\R^n)}~\leq~ \C_{p~\sigma~\Phi}~2^{\left({m\over n}\right)\ell}~\left\|f\right\|_{\L^p(\R^n)},\qquad {1\over p}~=~{1\over 2}(1-\vartheta).
\eeq
Recall $\vartheta=-{2m\over n-1}$. We have $\F^*_{\ell~\vartheta}=\F^*_\ell$ and 
${1\over 2}-{1\over p}={-m\over n-1}$. By duality, this implies
\bel{F_t L^p dual}
\begin{array}{lr}\ds
\left\|\F_\ell f\right\|_{\L^{p\over p-1}(\R^n)}~\leq~ \C_{p~\sigma~\Phi}~2^{\left({m\over n}\right)\ell}~\left\|f\right\|_{\L^{p\over p-1}(\R^n)}
\\\\ \ds
\hbox{for}\qquad{p-1\over p}-{1\over 2}~=~{-m\over n-1}~=~{1\over 2}-{1\over p}.
\end{array}
\eeq
By using (\ref{L^p Result vartheta}) and (\ref{F_t L^p dual}) and taking into account that $\sigma\in\S^{m}$ implies $\sigma\in\S^{m_1}$ for $m\leq m_1$,
we have
\bel{F_t L^p>}
\begin{array}{cc}\ds
\left\|\F_\ell f\right\|_{\L^p(\R^n)}~\leq~ \C_{p~\sigma~\Phi}~2^{\left({m\over n}\right)\ell}~\left\|f\right\|_{\L^p(\R^n)},\qquad\ell\ge0
\\\\ \ds
\hbox{whenever}\qquad \left|{1\over 2}-{1\over p}\right|~\leq~{-m\over n-1}.
\end{array}
\eeq
For $\ell\leq0$, we carry out same estimates in (\ref{Local est})-(\ref{F_t L^p dual}) except for (\ref{F_t a Result p 2}) and (\ref{F_t a Result 2 p' and F^*_t a Result p 2}) where we apply (\ref{F_t 2,p result<}) and (\ref{F_t p',2 result<}) instead of  (\ref{F_t 2,p result>}) and (\ref{F_t p',2 result>}) respectively. 

As a result, for $\sigma\in\S^m$, we have
\bel{F_t L^p<}
\begin{array}{cc}\ds
\left\|\F_\ell f\right\|_{\L^p(\R^n)}~\leq~ \C_{p~\sigma~\Phi}~2^{-\left({m\over n}\right)(n-1)\ell}~\left\|f\right\|_{\L^p(\R^n)},\qquad\ell\leq0
\\\\ \ds
\hbox{whenever}\qquad \left|{1\over 2}-{1\over p}\right|~\leq~{-m\over n-1}.
\end{array}
\eeq

\section{Majorization of Kernels}
\setcounter{equation}{0}
Let $\delta_\ell(\xi)$ be defined in (\ref{delta_t}) for every $\ell\in\Z$ and $\phi_j(\xi)$ be defined in (\ref{phi_j}) for every $j\in\Z$. Recall from (\ref{Partial}). We define
\bel{Omega_t,j}
\Omega_{\ell j}(x,y)~=~\int_{\R^n} e^{2\pi\i\left(\Phi(x,\xi)-y\cdot\xi\right)}\delta_\ell(\xi)\phi_j(\xi)\sigma(x,y,\xi)d\xi.
\eeq
Note that $\sigma(x,y,\xi)$ has a compact support in both $x$ and $y$. We abbreviate $\supp \sigma$ to be the support of $\sigma(x,y,\xi)$ in $x$ and $y$.

In order to prove (\ref{Comple est}), we write
\bel{F_t f rewrite j}
\begin{array}{lr}\ds
\int_{{^c}\Q_r(x_o)}\left| \Big(\F_\ell a\Big)(x)\right| dx~=~\int_{{^c}\Q_r(x_o)\cap\supp\sigma}\left|\int_{\R^n}a(y)\Omega_\ell(x,y)dy\right| dx
\\\\ \ds~~~~~~~~~~
~=~\int_{{^c}\Q_r(x_o)\cap\supp\sigma}\Bigg|\sum_{j\in\Z}\int_{\R^n} a(y)\Omega_{\ell j}(x,y)dy\Bigg|dx
\\\\ \ds~~~~~~~~~~
~\leq~\int_{{^c}\Q_r(x_o)\cap\supp\sigma}\left\{\int_{\R^n} |a(y)|\Bigg|\sum_{j\leq0}\Omega_{\ell j}(x,y)\Bigg|dy\right\}
+\sum_{j\ge0}\left|\int_{\R^n} a(y)\Omega_{\ell j}(x,y)dy\right|dx
\end{array}
\eeq
where $a$ is an $\H^1$-{\it atom} associated to the ball $B_r(x_o)$.

Observe that $\sum_{j\leq0}\phi_j(\xi)$ has a compact support  ( $|\xi|\leq2$ ) 
whereas  $|\supp\delta_\ell(\xi)|\leq\C2^{1-\ell}2^{(n-1)}$ for $\ell\ge0$ and $|\supp\delta_\ell(\xi)|\leq\C2^{1}2^{(n-1)\ell}$ for $\ell\leq0$. 
For $\sigma\in\S^0$, we have
\bel{Omega sum j<0}
\begin{array}{lr}\ds
\Bigg|\sum_{j\leq0}\Omega_{\ell j}(x,y)\Bigg|~=~\left|\int_{\R^n} e^{2\pi\i\left(\Phi(x,\xi)-y\cdot\xi\right)}\delta_\ell(\xi)\left\{\sum_{j\leq0}\phi_j(\xi)\right\}\sigma(x,y,\xi)d\xi\right|
\\\\ \ds~~~~~~~~~~~~~~~~~~~~~~~
~\leq~\int_{\R^n} \left|\delta_\ell(\xi)\right|\Bigg|\sum_{j\leq0}\phi_j(\xi)\Bigg|\left|\sigma(x,y,\xi)\right|d\xi
\\\\ \ds~~~~~~~~~~~~~~~~~~~~~~~
~\leq~\C \left\{\begin{array}{lr}\ds 2^{-\ell},\qquad \ell\ge0,
\\\\ \ds
2^{(n-1)\ell},\qquad \ell\leq0.\end{array}\right.
\end{array}
\eeq
By using (\ref{Omega sum j<0}), we have
\bel{Int Omega sum j<0}
\begin{array}{lr}\ds
\int_{{^c}\Q_r(x_o)\cap\supp\sigma}\left\{\int_{\R^n} |a(y)|\Bigg|\sum_{j\leq0}\Omega_{\ell j}(x,y)\Bigg|dy\right\}dx
\\\\ \ds
~\leq~
\left\{\begin{array}{lr}\ds
\C~\int_{\supp\sigma}\left\{2^{-\ell}\int_{\R^n} |a(y)|dy\right\}dx
~\leq~\C_\sigma~2^{-\ell}\left\|a\right\|_{\L^1(\R^n)},\qquad \ell\ge0,
\\\\ \ds
\C~\int_{\supp\sigma}\left\{2^{(n-1)\ell}\int_{\R^n} |a(y)|dy\right\}dx
~\leq~\C_\sigma~2^{(n-1)\ell}\left\|a\right\|_{\L^1(\R^n)},\qquad \ell\leq0.
\end{array}\right.
\end{array}
\eeq
\v

{\bf Principal Lemma~~}{\it  
Let $\Omega_{\ell j}(x,y)$ be defined in (\ref{Omega_t,j}) and $\sigma\in\S^{-{n-1\over 2}}$. For every $ j\ge0$, we have
\bel{Est1}
\int_{\R^n} \left|\Omega_{\ell j}(x,y)\right| dx~\leq~\C_{\sigma~\Phi}~\left\{\begin{array}{lr}\ds 
2^{-\left({n-1\over 2n}\right)\ell},\qquad \ell\ge0,
\\\\ \ds
2^{(n-1)\left({n-1\over 2n}\right)\ell},\qquad \ell\leq0,
\end{array}\right.
\eeq
\bel{Est2}
\int_{\R^n}\left|\Omega_{\ell j}(x,y)-\Omega_{\ell j}(x,x_o)\right|dx~\leq~\C_{\sigma~\Phi}~2^j|y-x_o|~~\left\{\begin{array}{lr}\ds 
2^{-\left({n-1\over 2n}\right)\ell},\qquad \ell\ge0,
\\\\ \ds
2^{(n-1)\left({n-1\over 2n}\right)\ell},\qquad \ell\leq0,
\end{array}\right.
\eeq
and
\bel{Est3}
y\in B_r(x_o),\qquad
\int_{{^c}\Q_r(x_o)}\left|\Omega_{\ell j}(x,y)\right|dx~\leq~\C_{\sigma~\Phi}~{2^{-j}\over r}~\left\{\begin{array}{lr}\ds 
2^{-\left({n-1\over 2n}\right)\ell},\qquad \ell\ge0,
\\\\ \ds
2^{(n-1)\left({n-1\over 2n}\right)\ell},\qquad \ell\leq0
\end{array}\right.
\eeq
whenever $2^{j}>r^{-1}$.}
\v
Let $j\ge0$. For $2^j\leq r^{-1}$,  we write
\bel{F_t cancella}
\int_{\R^n} a(y)\Omega_{\ell j}(x,y)dy~=~\int_{B_r(x_o)} a(y) \left(\Omega_{\ell j}(x,y)-\Omega_{\ell j}(x,x_o)\right)dy.
\eeq 
Note that 
$\int_{B_r(x_o)} a(y)dy=0$
whenever $a$ is an $\H^1$-atom.  

By using (\ref{Est2}), we have 
\bel{Norm Est1}
\begin{array}{lr}\ds
\int_{{^c}\Q_r(x_o)}\left|\int_{B_r(x_o)} a(y) \left(\Omega_{\ell j}(x,y)-\Omega_{\ell j}(x,x_o)\right)dy\right|dx
\\\\ \ds
~\leq~\int_{B_r(x_o)} |a(y)|\left\{\int_{\R^n} \left|\Omega_{\ell j}(x,y)-\Omega_{\ell j}(x,x_o)\right| dx\right\} dy
\\\\ \ds
~\leq~ \C_{\sigma~\Phi}~2^j|y-x_o|~2^{-\left({n-1\over 2n}\right)\ell}\|a\|_{\L^1(\R^n)}
\\\\ \ds
~\leq~\C_{\sigma~\Phi}~2^j r~2^{-\left({n-1\over 2n}\right)\ell}\|a\|_{\L^1(\R^n)},\qquad y\in B_r(x_o).
\end{array}
\eeq
By summing over all such $j$ s, we have 
\bel{Sum1}
\begin{array}{lr}\ds
\sum_{2^j\leq r^{-1}} \int_{{^c}\Q_r(x_o)}\left|\int_{\R^n} a(y)\Omega_{\ell j}(x,y)dy\right|dx
\\\\ \ds
~\leq~\C_{\sigma~\Phi}~\Bigg(\sum_{2^j\leq r^{-1}} 2^j \Bigg) ~ r~2^{-\left({n-1\over 2n}\right)\ell}\|a\|_{\L^1(\R^n)}\qquad \hbox{\small{by (\ref{F_t cancella})-(\ref{Norm Est1})}}
\\\\ \ds
~\leq~\C_{\sigma~\Phi}~2^{-\left({n-1\over 2n}\right)\ell}\|a\|_{\L^1(\R^n)}.
\end{array}
\eeq
For $2^j>r^{-1}$, by using (\ref{Est3}), we have
\bel{Norm Est2}
\begin{array}{lr}\ds
\int_{{^c}\Q_r(x_o)}\left|\int_{\R^n} a(y)\Omega_{\ell j}(x,y)dy\right|dx
\\\\ \ds
~\leq~\int_{B_r(x_o)} |a(y)|\left\{\int_{{^c}\Q_r(x_o)} \left|\Omega_{\ell j}(x,y)\right| dx\right\} dy
\\\\ \ds
~\leq~\C_{\sigma~\Phi}~{2^{-j}\over r}~2^{-\left({n-1\over 2n}\right)\ell}\left\|a\right\|_{\L^1(\R^n)}.
\end{array}
\eeq
By summing over all such $j$ s, we have
\bel{Sum2}
\begin{array}{lr}\ds
\sum_{2^j>r^{-1}} \int_{{^c}\Q_r(x_o)}\left|\int_{\R^n} a(y)\Omega_{\ell j}(x,y)dy\right|dx
\\\\ \ds
~\leq~\C_{\sigma~\Phi}~\Bigg(\sum_{2^j>r^{-1}} 2^{-j}\Bigg)~r^{-1}~2^{-\left({n-1\over 2n}\right)\ell}\left\|a\right\|_{\L^1(\R^n)}
\\\\ \ds
~\leq~\C_{\sigma~\Phi}~2^{-\left({n-1\over 2n}\right)\ell}\left\|a\right\|_{\L^1(\R^n)}.
\end{array}
\eeq
For $\ell\leq0$, we repeat all estimates in (\ref{Norm Est1})-(\ref{Sum2}) with  $2^{-\left({n-1\over 2n}\right)\ell}$ replaced by $2^{(n-1)\left({n-1\over 2n}\right)\ell}$.

On the other hand,  consider
\bel{F_t^* f rewrite j}
\begin{array}{lr}\ds
\int_{{^c}\Q^*_r(x_o)}\left| \Big(\F^*_\ell a\Big)(x)\right| dx~=~\int_{{^c}\Q^*_r(x_o)}\left|\int_{\R^n}a(y)\Omega^*_\ell(x,y)dy\right| dx
\end{array}
\eeq
where $\Omega^*_\ell(x,y)$ is defined in
 (\ref{F*_t singular}) and $\Q^*_r(x_o)$ is the {\it region of influence} associated to $\F^*_\ell$. 
 
Define
\bel{Omega^*_t,j}
\Omega^*_{\ell j}(x,y)~=~\int_{\R^n} e^{2\pi\i\left(x\cdot\xi-\Phi(y,\xi)\right)}\bar{\delta}_\ell(\xi)\phi_j(\xi)\bar{\sigma}(x,y,\xi)d\xi.
\eeq
Observe that same estimates hold in (\ref{F_t f rewrite j})-(\ref{Int Omega sum j<0}) for $\Omega^*_\ell(x,y)$ and $\Omega^*_{\ell j}(x,y), j\leq0$. Moreover,  {\bf Principal Lemma} is true for $\Omega^*_{\ell j}(x,y)$.

Let $\sigma\in\S^{-{n-1\over 2}}$. For every $j\ge0$, we have
\bel{Est1 dual}
\int_{\R^n} \left|\Omega^*_{\ell j}(x,y)\right| dx~\leq~\C_{\sigma~\Phi}~\left\{\begin{array}{lr}\ds 
2^{-\left({n-1\over 2n}\right)\ell},\qquad \ell\ge0,
\\\\ \ds
2^{(n-1)\left({n-1\over 2n}\right)\ell},\qquad \ell\leq0,
\end{array}\right.
\eeq
\bel{Est2 dual}
\int_{\R^n}\left|\Omega^*_{\ell j}(x,y)-\Omega^*_{\ell j}(x,x_o)\right|dx~\leq~\C_{\sigma~\Phi}~2^j|y-x_o|~~\left\{\begin{array}{lr}\ds 
2^{-\left({n-1\over 2n}\right)\ell},\qquad \ell\ge0,
\\\\ \ds
2^{(n-1)\left({n-1\over 2n}\right)\ell},\qquad \ell\leq0,
\end{array}\right.
\eeq
and
\bel{Est3 dual}
\begin{array}{lr}\ds
y\in B_r(x_o),\qquad
\int_{{^c}\Q^*_r(x_o)}\left|\Omega^*_{\ell j}(x,y)\right|dx~\leq~\C_{\sigma~\Phi}~{2^{-j}\over r}~\left\{\begin{array}{lr}\ds 
2^{-\left({n-1\over 2n}\right)\ell},\qquad \ell\ge0,
\\\\ \ds
2^{(n-1)\left({n-1\over 2n}\right)\ell},\qquad \ell\leq0
\end{array}\right.
\\ \ds
\hbox{whenever $2^{j}>r^{-1}$.}
\end{array}
\eeq

We obtain (\ref{H^1 est *})  by carrying out same estimates in (\ref{F_t cancella})-(\ref{Sum2}) with $\Omega_{\ell j}(x,y)$ replaced by $\Omega^*_{\ell j}(x,y)$ and applying (\ref{Est1 dual})-(\ref{Est3 dual}) instead.

\section{A second dyadic decomposition}
\setcounter{equation}{0}
By  {\it localization principal} of oscillatory integrals given in chapter VIII of \cite{Stein},
the kernel $\Omega_{\ell j}(x,y)$ defined in (\ref{Omega_t,j}) has singularity appeared at 
\bel{localization}
\nabla_\xi\Big(\Phi(x,\xi)-y\cdot\xi\Big)~=~0.
\eeq
For every $x\in\R^n$, we consider the variety 
\bel{locus}
\Sigma_x~=~\Big\{~ y\in\R^n~\colon~y=\nabla_\xi\Phi(x,\xi)~~\hbox{for some}~\xi ~\Big\}
\eeq
which is the locus of the singularity of $y\mt\Omega_\ell(x,y)$. 

Note that $\nabla_\xi\Phi(x,\xi)$ is homogeneous of degree zero in $\xi$. The projection of $\Sigma_x$ on the unit sphere  $\mathbb{S}^{n-1}\subset\R^n$  has dimension at most equal to $n-1$. 

Let $j\ge0$ be fixed.   
We construct a set of points,  denoted by $\left\{\xi^\nu_j\right\}_\nu$  , that are almost equally distributed on  $\mathds{S}^{n-1}$ with grid length equal to $\mathfrak{B} 2^{-j/2}$ for ${1/2}\leq\mathfrak{B}\leq2$. \begin{remark}
For every $\xi\in\R^n$, there exists a $\xi^\nu_j\in\left\{\xi^\nu_j\right\}_\nu$ such that $\left|{\xi\over|\xi|}-\xi^\nu_j\right|\leq2^{-j/2}$.
\end{remark}

Consider
\bel{Intersection k}
\begin{array}{cc}\ds
\mathbb{S}^{n-2}_k~\doteq~\mathbb{S}^{n-1}\cap\left\{(\tau,\lambda)\in\R\times\R^{n-1}~\colon~\tau=k~2^{-j/2-\kappa}\right\},
\\\\ \ds
\kappa=\left\{\begin{array}{lr}\ds {1/ 2}~~~~\hbox{if $j$ is odd}
\\\\ \ds
~~1~~~~~\hbox{if $j$ is even}
\end{array}\right.,
\qquad   k=-2^{j/2+\kappa},\ldots,0,\ldots,2^{j/2+\kappa}.
 \end{array}
\eeq
\v
Observe that for $k=2^{j/2+\kappa}$ or $-2^{j/2+\kappa}$, the intersection $\mathbb{S}^{n-2}_k$ is a point $(\tau,\lambda)=(1,0)$ or $(-1,0)$. On the other hand, 
 $\mathbb{S}^{n-2}_0=\mathbb{S}^{n-2}$ is the unit sphere in the $(n-1)$-dimensional $\lambda$-space.

For each $k$, we choose a set of points that are equally distributed on $\mathbb{S}^{n-2}_k$ with grid length equal to $\mathfrak{B} 2^{-j/2}$ for some $1/2\leq\mathfrak{B}\leq1$. Define $\left\{\xi^\nu_j\right\}_\nu$ to be
the union of these sets, for all $k=-2^{j/2+\kappa},\ldots,0,\ldots,2^{j/2+\kappa}$.
It is clear that there are at most $\C2^{j\left({n-1\over2}\right)}$ elements in  $\left\{\xi^{\nu}_j\right\}_{\nu}$. 
Moreover, the following two conditions are satisfied.

{\bf ( 1 )} A subset of points in $\left\{\xi^\nu_j\right\}_\nu$ are equally distributed on $\mathbb{S}^{n-2}$ in the ($n-1$)-dimensional $\lambda$-space,  with grid length equal to $\mathfrak{B}2^{-j/2}$ for $1/2\leq\mathfrak{B}\leq1$.

{\bf ( 2 )}  $(\tau,\lambda)=(\pm1,0)\in\R\times\R^{n-1}$ belong to the collection $\left\{\xi^{\nu}_j\right\}_{\nu}$.

Let $\delta_\ell(\xi)$ be defined in (\ref{delta_t}) which is supported on the dyadic cone $\Lambda_\ell$ given in (\ref{Cone}). Consider $\xi=(\tau,\lambda)\in\Lambda_\ell\cap\mathbb{S}^{n-1}$. 
For $\ell>{j/2}+3$, we have $|\tau|<2^{1-\ell+1}<2^{-j/2-1}$ and for
 $\ell<-j/2-3$, we have $|\lambda|<2^{1-\ell+1}<2^{-j/2-1}$.  
 We can verify the following:
 \begin{remark}
 For every $\xi\in\Lambda_\ell, ~\ell>j/2+3$, there exists a $\xi^\nu_j\in\mathbb{S}^{n-2}$ in the $\lambda$-space such that  $\left|{\xi\over|\xi|}-\xi^\nu_j\right|\leq2^{-j/2}$.
 
 For every $\xi\in\Lambda_\ell, ~\ell<-j/2-3$, there is a $\xi^\nu_j=(\tau,\lambda)=(\pm1,0)$ such that  $\left|{\xi\over|\xi|}-\xi^\nu_j\right|\leq2^{-j/2}$. 
 \end{remark}
 
 Let $\Gamma^{\nu}_{j}$ denote the cone whose central direction is $\xi^{\nu}_{j}$, such that 
\bel{Gamma}
\Gamma_{j}^{\nu}~=~\left\{\xi\in\R^n~\colon~\left| {\xi\over |\xi|}-\xi^{\nu}_{j}\right|~\leq~2\cdot2^{-j/2}\right\}.
\eeq
Let $\varphi$ be defined in (\ref{varphi funda}). Consider
\bel{phi^v_j}
\varphi^\nu_j(\xi)~=~\varphi\left(2^{j/2}\left({\xi\over |\xi|}-\xi^{\nu}_{j}\right)\right)
\eeq
which is supported on $\Gamma^\nu_j$.

Let $\phi_j(\xi)$ be defined in (\ref{phi_j}) where $2^{j-1}\leq|\xi|<2^{j+1}$ for $\xi$ in the support of $\phi_j(\xi)$. From (\ref{Gamma}) we have
\bel{Gamma norm}
\left|\Gamma_j^\nu\cap\left\{2^{j-1}\leq|\xi|<2^{j+1}\right\}\right|~\leq~\C~2^j2^{j\left({n-1\over 2}\right)}.
\eeq
On the other hand, by definition of $\Lambda_\ell$ in (\ref{Cone}), we have
\bel{Lambda norm><}
\begin{array}{cc}\ds
\left|\Lambda_\ell\cap\left\{2^{j-1}\leq|\xi|<2^{j+1}\right\}\right|~\leq~\C~ 2^{j-\ell}2^{j\left({n-1}\right)},\qquad \ell\ge0,
\\\\ \ds
\left|\Lambda_\ell\cap\left\{2^{j-1}\leq|\xi|<2^{j+1}\right\}\right|~\leq~\C~ 2^{j}2^{(j+\ell)\left({n-1}\right)},\qquad \ell\leq0.
\end{array}
\eeq

We begin to construct  a smooth partition of unity, associated  $\Gamma^\nu_j$  in (\ref{Gamma}),  respectively for 
 \bel{Cases}
\begin{array}{cc}
\hbox{\bf Case One:}\qquad  -j/2-3~\leq~\ell~\leq~ j/2+3,
\\\\ \ds
\hbox{\bf Case Two:}\qquad  \ell~>~j/2+3,\qquad \hbox{\bf Case Three:}\qquad \ell~<-j/2-3.
\end{array}
\eeq
\v

{\bf Case One:} Let $-j/2-3\leq\ell\leq{j/2}+3$. Define
\bel{chi^v_j}
\chi^\nu_{\ell j}(\xi)~=~\varphi^\nu_j(\xi)\Big/\sum_{\nu}\varphi^\nu_j(\xi),\qquad \xi\in\mathbb{C}
\eeq
where the summation is taking over all elements in $\left\{\xi^\nu_j\right\}_\nu$. 
By {\bf Remark 6.1}, we have
\bel{chi prop}
\sum_{\nu} \chi^{\nu}_{j}\left(\xi\right)=1,\qquad\xi\in\mathbb{C}.
\eeq
Suppose $\Gamma^{\nu}_j\cap\Lambda_\ell\neq\emptyset$. We  clearly have
\bel{intersection of cone norm}
\left|\Gamma^\nu_j\cap\Lambda_\ell\cap\left\{2^{j-1}\leq|\xi|<2^{j+1}\right\}\right|~\leq~\C~2^j2^{j\left({n-1\over 2}\right)}\qquad\hbox{by (\ref{Gamma norm}).}
\eeq
Moreover, by definition of $\Lambda_\ell$ in (\ref{Cone}), the area
\bel{area}
\left|\Lambda_\ell\cap\mathbb{S}^{n-1}\right|~\leq~\left\{\begin{array}{lr}\ds \C ~2^{-\ell},\qquad~~~~~ \ell\ge0,
\\\\ \ds
\C~2^{\ell(n-1)},\qquad \ell\leq0.
\end{array}\right.
\eeq
Note that $\left\{\xi^{\nu}_j\right\}_\nu$ are  almost equally distributed on $\mathbb{S}^{n-1}$ with grid length  $\mathfrak{B} 2^{-j/2}$ for ${1\over 2}\leq\mathfrak{B}\leq2$.
There are at most 
\bel{number elements case3}
\left\{\begin{array}{lr}\ds
\C~2^{j\left({n-1\over 2}\right)}2^{-\ell},\qquad ~~~~~0\leq\ell\leq{j/ 2}+3,
\\\\ \ds
\C~2^{j\left({n-1\over 2}\right)}2^{\ell(n-1)},\qquad -{j/ 2}-3\leq\ell\leq0
\end{array}\right.
\eeq
elements in $\left\{\xi^{\nu}_j\right\}_\nu$ such that $\Lambda_\ell\cap\Gamma^\nu_j$ is nonempty.

\v

{\bf Case Two:} Let $\ell>{j/2}+3$. We have $\Lambda_\ell$ intersect with $\Gamma^\nu_j$ whose central direction     is some $\xi^\nu_j\in\mathbb{S}^{n-2}$ in the  $\lambda$-space. Define
\bel{chi^v_j>}
\chi^\nu_{\ell j}(\xi)~=~ \varphi^\nu_j(\xi)~\Big/\sum_{\nu~\colon~\xi^\nu_j\in\mathbb{S}^{n-2}} \varphi^\nu_j(\xi),\qquad  \xi\in\Lambda_\ell.
\eeq
By {\bf Remark 6.2},  we have
\bel{chi prop >}
\sum_{\nu~\colon~\xi^\nu_j\in\mathbb{S}^{n-2}} \chi^{\nu}_{\ell j}\left(\xi\right)~=~1,\qquad \xi\in\Lambda_\ell.
\eeq
Suppose $\Gamma^{\nu}_j\cap\Lambda_\ell\neq\emptyset$. By putting together (\ref{Gamma norm}) and (\ref{Lambda norm><}),  we have
\bel{intersection of cone norm>}
\left|\Gamma^\nu_j\cap\Lambda_\ell\cap\left\{2^{j-1}\leq|\xi|<2^{j+1}\right\}\right|~\leq~\C~2^{j-\ell}2^j2^{j\left({n-2\over 2}\right)}.
\eeq
There are at most $\C ~2^{\left({n-2\over 2}\right)j}$ such elements in $\left\{\xi^{\nu}_j\right\}_\nu$ equally distributed on $\mathbb{S}^{n-2}$.

\v

{\bf Case Three:} Let $\ell<-j/2-3$. We have $\Lambda_\ell$ intersect  with $\Gamma^\nu_j$ whose central direction  is $\xi^\nu_j=(\tau,\lambda)=(\pm1,0)$. Define
\bel{chi^v_j<}
\chi^\nu_{\ell j}(\xi)~=~\varphi^\nu_j(\xi)~\Big/\sum_{\nu~\colon~ \xi^\nu_j=(\pm1,0)} \varphi^\nu_j(\xi),\qquad 
\xi\in\Lambda_\ell.
\eeq
By  {\bf Remark 6.2},  we have
 \bel{chi prop <}
\sum_{\nu~\colon~\xi^\nu_j=(\pm1,0)} \chi^{\nu}_{\ell j}\left(\xi\right)~=~1,\qquad \xi\in\Lambda_\ell.
\eeq
By putting together (\ref{Gamma norm}) and (\ref{Lambda norm><}),  we  have
\bel{intersection of cone norm<}
\left|\Gamma^\nu_j\cap\Lambda_\ell\cap\left\{2^{j-1}\leq|\xi|<2^{j+1}\right\}\right|~\leq~\C~2^j2^{(j+\ell)\left({n-1}\right)}.
\eeq

\begin{figure}[h]
\centering
\includegraphics[scale=0.40]{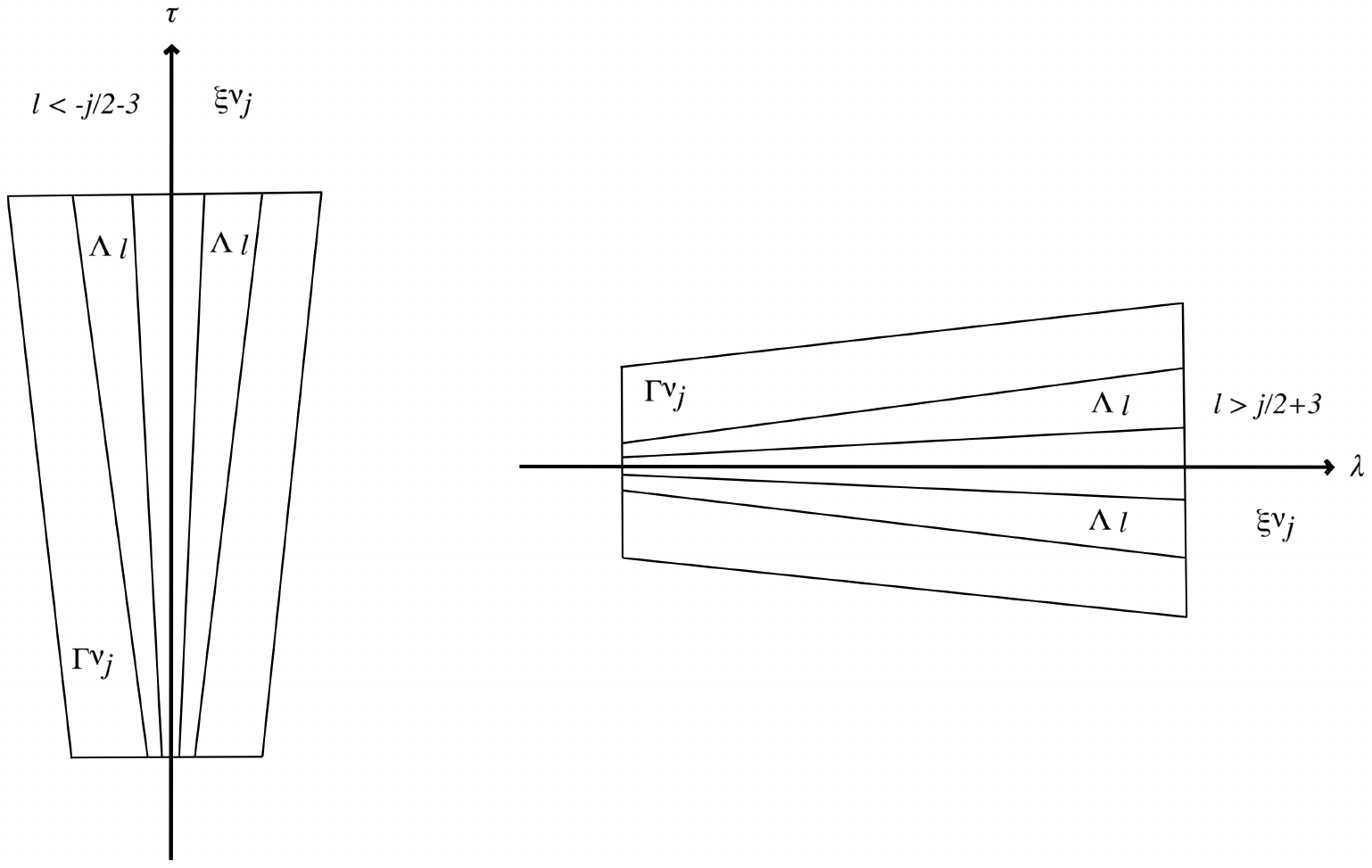}
\caption{the  enclosed regions are examples of 
$\Gamma^\nu_j\cap\Lambda_\ell\cap\left\{2^{j-1}\leq|\xi|<2^{j+1}\right\}$ for $\ell<-j/2-3$ where $\xi^\nu_j=(\tau,\lambda)=(1,0)$ and for $\ell> j/2+3$ where $\xi^\nu_j\in\mathbb{S}^{n-2}$ in the $\lambda$-space.}
\end{figure}
Next, for computational purposes, we introduce a linear isometry. For every $\nu$, consider  $\xi=\L_\nu \eta$ where $\L_\nu$ is an $n\times n$-matrix of rotations with $\det\L_\nu=1$. 

In particular,  we require that the $\imath$-th coordinate of $\eta$ is in the same direction of $\xi^\nu_j$ for some $\imath\in\{1,2,\ldots,n\}$.  Denote 
\bel{eta^v_j}
\eta^\nu_j~\doteq~\left({\eta_\imath\over|\eta_\imath|},0\right)~\in~\R\times\R^{n-1}.
\eeq
We thus have
 $\xi^\nu_j=\L_\nu \eta^\nu_j$.  
 
Without lose of the generality, we fix $\tau=\xi_1\in\R$ and $\lambda=\xi_1^\dagger\in\R^{n-1}$.
For condition {\bf ( 1 )} where $\xi^\nu_j\in\mathbb{S}^{n-2}$ in the $\lambda$-space,  we require
\bel{L_v decomposable}
\L_\nu~=~\left[\begin{array}{ccc}
1 \qquad
\\ \ds
\qquad \L'_\nu
\end{array}\right],\qquad \det\L'_\nu~=~1\qquad \hbox{and}\qquad \eta_\imath~=~\eta_n
\eeq
where $\L'_\nu$ is an $(n-1)\times(n-1)$-matrix. 

For condition {\bf ( 2 )} where $\xi^\nu_j=(\tau,\lambda)=(\pm1,0)$, we choose
$\L_\nu$ to be the identity matrix so that 
\bel{eta_1}
\eta_\imath~=~\eta_1~=~\xi_1~=~\tau.
\eeq 
Let  $\varphi^\nu_j(\xi)$ be defined in (\ref{phi^v_j}) and $\chi^\nu_{\ell j}(\xi)$ be defined respectively in (\ref{chi^v_j}), (\ref{chi^v_j>}) and (\ref{chi^v_j<}). Observe that
\bel{chi prop est}
\begin{array}{cc}\ds
\left|\left({\p\over\p\eta}\right)^\alpha \chi_{\ell j}^{\nu}\left(\L_\nu\eta\right)\right|~\leq~\C_{\alpha} ~2^{|\alphaup| \left({1\over 2}\right)j}|\eta|^{-|\alphaup|}
\end{array}
\eeq
for every multi-index $\alpha$. 

Let $r=|\xi|=\left|\L_\nu\eta\right|=|\eta|$.   For every $\L_\nu\eta=\xi\in\Gamma^\nu_j$, the angle between $\eta$ and $\eta_\imath$ is bounded by $ \arcsin(2\cdot2^{-j/2})$.
From direct computation, we have
\bel{radial derivative}
{\p\over \p\eta_\imath}~=~\left({\p r\over \p \eta_\imath}\right){\p\over \p r}+\O\left(2^{-j/2}\right)\cdot\nabla_{\eta_\imath^\dagger}.
\eeq
Note that $\chi^\nu_{\ell j}(\xi)=\chi^\nu_{\ell j}\left(\L_\nu \eta\right)$  is homogeneous of degree zero in $\eta$. Hence that $\p_r\chi^\nu_{\ell j}\equiv0$.
Together with the estimate in (\ref{chi prop est}), we have
\bel{chi prop est split}
\begin{array}{cc}\ds
\left|\left({\p\over \p\eta_\imath}\right)^\alpha \chi_{\ell j}^{\nu}\left(\L_\nu\eta\right)\right|~\leq~\C_{\alpha} ~|\eta|^{-|\alpha|},
\qquad
\left|\left({\p\over\p\eta_\imath^\dagger}\right)^\beta \chi_{\ell j}^{\nu}\left(\L_\nu\eta\right)\right|~\leq~\C_{\beta} ~2^{|\beta| \left({1\over2}\right)j}|\eta|^{-|\beta|}
\end{array}
\eeq
for every multi-indices $\alpha, \beta$.
\v
Now, we  explicitly construct our {\it region of influence} $\mathfrak{Q}_r(x_o)$ discussed earlier in section 4.

{\bf Case One:} Let $-j/2-3\leq\ell\leq j/2+3$. Consider 
\bel{rectangle R}
\begin{array}{rl}\ds
R^\nu_j(x_o)~=~\left\{x\in\R^n\colon
\left|\left(\L_\nu^Tx_o-\nabla_\eta\Phi\left(x,\L_\nu \eta^\nu_j\right)\right)_\imath\right|\leq4\cdot2^{-j},~
\left|\left(\L_\nu^Tx_o-\nabla_\eta\Phi\left(x,\L_\nu\eta^\nu_j\right)\right)_\imath^\dagger\right|\leq 4\cdot2^{-j/2}
\right\}
\end{array}
\eeq
for every $\xi^\nu_j$ in the collection $\left\{\xi^\nu_j\right\}_\nu$.  ( $\L_\nu^T$ is the transpose of $\L_\nu$. ) 

The set $\mathfrak{Q}_r(x_o)$ is defined by
\bel{Q_r}
\mathfrak{Q}_r(x_o)~=~\bigcup_{2^{-j}\leq r}~\Bigg(~\bigcup_{\nu} R_{j}^{\nu}(x_o)~\Bigg).
\eeq
We have 
\bel{Q_r Est exact}
\begin{array}{lr}\ds
 \left|\Q_r(x_o)\right|~\leq~\sum_{2^{-j}\leq r}\sum_{\nu}\left|R^{\nu}_{j}(x_o)\right|
 \\\\ \ds ~~~~~~~~~~~~~
 ~\leq~\C~\sum_{2^{-j}\leq r}\sum_{\nu} 2^{-j\left({n-1\over 2}\right)}2^{-j}
 ~=~\C~\sum_{2^{-j}\leq r}\sum_{\nu} 2^{-j\left({n+1\over 2}\right)}
 \\\\ \ds ~~~~~~~~~~~~~
 ~\leq~\C~\sum_{2^{-j}\leq r} 2^{j\left({n-1\over 2}\right)}2^{-j\left({n+1\over 2}\right)} ~=~\C~\sum_{2^{-j}\leq r} 2^{-j}
 \\\\ \ds ~~~~~~~~~~~~~
 ~\leq~\C~r
  \end{array}
\eeq
satisfying the estimate in (\ref{Q norm}).   

Recall that $\sigma(x,y,\xi)$ has a compact support in both $x$ and $y$.

{\bf Case Two:} Let $\ell>j/2+3$. Consider 
\bel{rectangle R>}
\begin{array}{rl}\ds
R^\nu_j(x_o)=
\left\{x\in\supp\sigma~\colon
\left|\left(\L_\nu^Tx_o-\nabla_\eta\Phi\left(x,\L_\nu\eta_j^\nu\right)\right)_n\right|\leq4\cdot2^{-j},~
\left|\left(\L_\nu^Tx_o-\nabla_\eta\Phi\left(x,\L_\nu\eta_j^\nu\right)\right)_1^\dagger\right|\leq 4\cdot2^{-j/2}
\right\}
\end{array}
\eeq
for every  $\xi^\nu_j\in\mathbb{S}^{n-2}$ in the $(n-1)$-dimensional $\lambda$-space. 

Note that $\eta_\imath=\eta_n$ and the first coordinate
$\left(\L_\nu^Tx_o-\nabla_\eta\Phi\left(x,\L_\nu\eta_j^\nu\right)\right)_1\in R^\nu_j(x_o)
$
given in (\ref{rectangle R>}) has no other restriction except for $x\in\supp \sigma$.

The set $\mathfrak{Q}_r(x_o)$ is defined by
\bel{Q_r>}
\mathfrak{Q}_r(x_o)~=~\bigcup_{2^{-j}\leq r}~\Bigg(~\bigcup_{\nu~\colon~\xi^\nu_j\in\mathbb{S}^{n-2}} R_{j}^{\nu}(x_o)~\Bigg).
\eeq
For $\xi^\nu_j\in\mathbb{S}^{n-2}$, there are at most $\C2^{j\left({n-2\over 2}\right)}$ such elements in $\left\{\xi^\nu_j\right\}_\nu$. We have
\bel{Q_r Est exact>}
\begin{array}{lr}\ds
 \left|\Q_r(x_o)\right|~\leq~\sum_{2^{-j}\leq r}~~\sum_{\nu~\colon~\xi^\nu_j\in\mathbb{S}^{n-2}}\left|R^{\nu}_{j}(x_o)\right|
 \\\\ \ds ~~~~~~~~~~~~~
 ~\leq~\C_\sigma~\sum_{2^{-j}\leq r}~~\sum_{\nu~\colon~\xi^\nu_j\in\mathbb{S}^{n-2}} 2^{-j\left({n-2\over 2}\right)}2^{-j}
 ~=~\C_\sigma~\sum_{2^{-j}\leq r}~~\sum_{\nu~\colon~\xi^\nu_j\in\mathbb{S}^{n-2}} 2^{-j\left({n\over 2}\right)}
 \\\\ \ds ~~~~~~~~~~~~~
 ~\leq~\C_\sigma~\sum_{2^{-j}\leq r} 2^{j\left({n-2\over 2}\right)}2^{-j\left({n\over 2}\right)} ~=~\C_\sigma~\sum_{2^{-j}\leq r} 2^{-j}
 \\\\ \ds ~~~~~~~~~~~~~
 ~\leq~\C_\sigma~r
  \end{array}
\eeq
also satisfying the estimate in (\ref{Q norm}).

{\bf Case Three:} Let $\ell<-j/2-3$. Consider
\bel{rectangle R<}
\begin{array}{rl}\ds
R^\nu_j(x_o)~=~\left\{x\in\supp\sigma~\colon~
\left|\left(\L_\nu^Tx_o-\nabla_\eta\Phi\left(x,\L_\nu\eta_j^\nu\right)\right)_1\right|~\leq~4\cdot2^{-j}
\right\}
\end{array}
\eeq
for which $\eta^\nu_j=\xi^\nu_j=(\pm1,0)$ where $\eta_\imath=\eta_1$. Observe  that  
$\left(\L_\nu^Tx_o-\nabla_\eta\Phi\left(x,\L_\nu\eta_j^\nu\right)\right)_1^\dagger\in R^\nu_j(x_o)$
given  in (\ref{rectangle R<}) has no other restriction except for $x\in\supp \sigma$.

The set $\Q_r(x_o)$ is then defined by
\bel{Q_r<}
\Q_r(x_o)~=~\bigcup_{2^{-j}\leq r}~\Bigg(~\bigcup_{\nu~\colon~\xi^\nu_j=(\pm1,0)} R_{j}^{\nu}(x_o)~\Bigg).
\eeq
It is clear that
\bel{Q_r Est exact<}
\begin{array}{lr}\ds
 \left|\Q_r(x_o)\right|~\leq~\sum_{2^{-j}\leq r}~~\sum_{\nu~\colon~\xi^\nu_j=(\pm1,0)}\left|R^{\nu}_{j}(x_o)\right|
 \\\\ \ds ~~~~~~~~~~~~~
 ~\leq~\C_\sigma~\sum_{2^{-j}\leq r}~~\sum_{\nu~\colon~\xi^\nu_j=(\pm1,0)} 2^{-j}
  \\\\ \ds ~~~~~~~~~~~~~
 ~\leq~\C_\sigma~\sum_{2^{-j}\leq r} 2^{-j}
 \\\\ \ds ~~~~~~~~~~~~~
 ~\leq~\C_\sigma~r.
  \end{array}
\eeq

On the other hand, for the adjoint operator $\F^*_\ell$, we define the associated $\mathfrak{Q}^*_r(x_o)$  as follows.

{\bf Case One:} Let $-j/2-3\leq\ell\leq j/2+3$. Consider 
\bel{rectangle R*}
\begin{array}{rl}\ds
{^*}R^\nu_j(x_o)~=~\left\{x\in\R^n\colon
\left|\left(\L_\nu^Tx-\nabla_\eta\Phi\left(x_o,\L_\nu \eta^\nu_j\right)\right)_\imath\right|\leq4\cdot2^{-j},~
\left|\left(\L_\nu^Tx-\nabla_\eta\Phi\left(x_o,\L_\nu\eta^\nu_j\right)\right)_\imath^\dagger\right|\leq 4\cdot2^{-j/2}
\right\}
\end{array}
\eeq
for every $\xi^\nu_j$ in the collection $\left\{\xi^\nu_j\right\}_\nu$. 

The set $\mathfrak{Q}_r^*(x_o)$ is defined by
\bel{Q_r^*}
\mathfrak{Q}_r^*(x_o)~=~\bigcup_{2^{-j}\leq r}~\Bigg(~\bigcup_{\nu} {^*}R_{j}^{\nu}(x_o)~\Bigg).
\eeq

{\bf Case Two:}  Let $\ell>j/2+3$. Consider 
\bel{rectangle R*>}
\begin{array}{rl}\ds
{^*}R^\nu_j(x_o)=
\left\{x\in\supp\sigma~\colon
\left|\left(\L_\nu^Tx-\nabla_\eta\Phi\left(x_o,\L_\nu\eta_j^\nu\right)\right)_n\right|\leq4\cdot2^{-j},~
\left|\left(\L_\nu^Tx-\nabla_\eta\Phi\left(x_o,\L_\nu\eta_j^\nu\right)\right)_1^\dagger\right|\leq 4\cdot2^{-j/2}
\right\}
\end{array}
\eeq
for every  $\xi^\nu_j\in\mathbb{S}^{n-2}$ in the $(n-1)$-dimensional $\lambda$-space. 

Note that $\eta_\imath=\eta_n$ and
$\left(\L_\nu^Tx-\nabla_\eta\Phi\left(x_o,\L_\nu\eta_j^\nu\right)\right)_1\in  {^*}R^\nu_j(x_o)$
given  in (\ref{rectangle R*>}) has no other restriction except for $x\in\supp \sigma$.

The set $\mathfrak{Q}_r^*(x_o)$ is defined by
\bel{Q_r^*>}
\mathfrak{Q}_r^*(x_o)~=~\bigcup_{2^{-j}\leq r}~\Bigg(~\bigcup_{\nu~\colon~\xi^\nu_j\in\mathbb{S}^{n-1}} {^*}R_{j}^{\nu}(x_o)~\Bigg).
\eeq

{\bf Case Three:} Let $\ell<-j/2-3$. Consider
\bel{rectangle R*<}
\begin{array}{rl}\ds
{^*}R^\nu_j(x_o)~=~\left\{x\in\supp\sigma~\colon~
\left|\left(\L_\nu^Tx-\nabla_\eta\Phi\left(x_o,\L_\nu\eta_j^\nu\right)\right)_1\right|~\leq~4\cdot2^{-j}
\right\}
\end{array}
\eeq
for which $\eta^\nu_j=\xi^\nu_j=(\pm1,0)$ where $\eta_\imath=\eta_1$.

Observe  that  
$\left(\L_\nu^Tx-\nabla_\eta\Phi\left(x_o,\L_\nu\eta_j^\nu\right)\right)_1^\dagger\in  {^*}R^\nu_j(x_o)$
given  in (\ref{rectangle R*<}) has no other restriction except for $x\in\supp \sigma$.

The set $\Q_r^*(x_o)$ is then defined by
\bel{Q_r^*<}
\Q_r^*(x_o)~=~\bigcup_{2^{-j}\leq r}~\Bigg(~\bigcup_{\nu~\colon~\xi^\nu_j=(\pm1,0)} {^*}R_{j}^{\nu}(x_o)~\Bigg).
\eeq
It is easy to verify that  same estimates  in (\ref{Q_r Est exact}), (\ref{Q_r Est exact>}) and (\ref{Q_r Est exact<}) hold for $\Q^*_r(x_o)$ defined respectively in (\ref{Q_r^*}), (\ref{Q_r^*>}) and (\ref{Q_r^*<}). Hence that we find (\ref{Q* norm}) as desired.

\section{Proof of Principal Lemma}
\setcounter{equation}{0}
Let $\Omega_{\ell j}(x,y)$ be defined in (\ref{Omega_t,j}). 
For every $j\ge0$, we aim to show  
\bel{Est1*}
 \int_{\R^n} \left|\Omega_{\ell j}(x,y)\right| dx~\leq~\C_{\sigma~\Phi}~\left\{\begin{array}{lr}\ds 2^{-\left({1\over 2}\right)\ell},\qquad ~~\ell\ge0,
 \\\\ \ds
 2^{\left({n-1\over 2}\right)\ell},~~~~~~~\ell\leq0;
 \end{array}\right. 
 \eeq
 
\bel{Est2*}
\int_{\R^n} \left|\Omega_{\ell j}(x,y)-\Omega_{\ell j}(x,x_o)\right| dx~\leq~\C_{\sigma~\Phi}~2^j|y-x_o|
 \left\{\begin{array}{lr}\ds 2^{-\left({1\over 2}\right)\ell},\qquad ~~\ell\ge0,
 \\\\ \ds
 2^{\left({n-1\over 2}\right)\ell},~~~~~~~ \ell\leq0
 \end{array}\right. 
 \eeq
 
 and
\bel{Est3*}
\int_{{^c}\mathfrak{Q}_r(x_o)} \left|\Omega_{\ell j}(x,y)\right| dx~\leq~\C_{\sigma~\Phi}~{2^{-j}\over r}
 \left\{\begin{array}{lr}\ds 2^{-\left({1\over 2}\right)\ell},\qquad ~~~~~\ell\ge0,
 \\\\ \ds
 2^{\left({n-1\over 2}\right)\ell},\qquad~~ \ell\leq0
 \end{array}\right. 
 \eeq
for $y\in B_r(x_o)$ whenever $2^j\ge r^{-1}$.
 
Note that (\ref{Est1*})-(\ref{Est3*}) implies (\ref{Est1})-(\ref{Est3}) respectively since $(n-1)/2n\leq 1/2$ for every $n\ge2$.

Let $\chi^\nu_{\ell j}(\xi)$ be defined respectively in (\ref{chi^v_j}), (\ref{chi^v_j>}) and (\ref{chi^v_j<}). We  write
\bel{Omega_lj Sum}
\Omega_{\ell j}(x,y)~=~\sum_\nu \Omega_{\ell j}^\nu(x,y),\qquad -j/2-3\leq\ell\leq j/2+3;
\eeq
\bel{Omega_lj Sum >}
\Omega_{\ell j}(x,y)~=~\sum_{\nu~\colon~\xi^\nu_j\in\mathbb{S}^{n-2}}\Omega_{\ell j}^\nu(x,y),\qquad \ell>{j/2}+3;
\eeq
\bel{Omega_lj Sum <}
\Omega_{\ell j}(x,y)~=~\sum_{\nu~\colon~ \xi^\nu_j=(\pm1,0)} \Omega_{\ell j}^\nu(x,y),\qquad \ell<-j/2-3
\eeq
where
\bel{Omega^v_lj}
\Omega^\nu_{\ell j}(x,y)~=~\int_{\R^n} e^{2\pi\i\left(\Phi(x,\xi)-y\cdot\xi\right)}\chi^\nu_{\ell j}(\xi)\delta_\ell(\xi)\phi_j(\xi)\sigma(x,y,\xi)d\xi.
\eeq

Recall from the previous section. $\L_\nu$ is a $n\times n$-matrix of rotations. We have $\xi=\L_\nu\eta$ and $\det\L_\nu=1$. In particular, $\xi^\nu_j=\L_\nu\eta^\nu_j$ for which $\ds\eta^\nu_j=\Big(\eta_\imath/|\eta_\imath|,0\Big)\in\R\times\R^{n-1}$ as (\ref{eta^v_j}). Consider
\bel{Phi split}
\begin{array}{cc}\ds
\Phi(x,\L_\nu\eta)-y\cdot \L_\nu\eta
~=~\left(\nabla_\eta\Phi\left(x,\L_\nu\eta_j^\nu\right)-\L_\nu^T y\right)\cdot\eta~+~\Psi(x,\eta),
\\\\ \ds
\Psi(x,\eta)~=~\Phi(x,\L_\nu\eta)-\nabla_\eta\Phi\left(x,\L_\nu\eta_j^\nu\right)\cdot\eta.
\end{array}
\eeq
From (\ref{Omega_lj Sum}), we rewrite
\bel{Omega rewrite}
\begin{array}{cc}\ds
\Omega^\nu_{\ell j}(x,y)~=~\int_{\R^n} e^{2\pi\i\left(\nabla_\eta\Phi\left(x,\L_\nu\eta_j^\nu\right)-\L_\nu^T y\right)\cdot\eta}\Theta^\nu_{\ell j}(x,y,\eta)d\eta,
\\\\ \ds
\Theta^\nu_{\ell j}(x,y,\eta)~=~e^{2\pi\i\Psi(x,\eta)}\chi^\nu_{\ell j}(\L_\nu\eta)\delta_\ell(\L_\nu\eta)\phi_j(\L_\nu\eta)\sigma(x,y,\L_\nu\eta).
\end{array}
\eeq
It is clear that for $\sigma\in\S^{-{n-1\over 2}}$, we have
\bel{Theta norm}
\left|\Theta^\nu_{\ell j}(x,y,\eta)\right|~\leq~\C~\left({1\over 1+|\eta|}\right)^{n-1\over 2}~\leq~\C~2^{-j\left({n-1\over 2}\right)}.
\eeq
Moreover, $\Theta^\nu_{\ell j}(x,y,\eta)$ has a compact support in both $x$ and $y$. ( $\supp \sigma$ )

Let $\Gamma^\nu_j$ be defined in (\ref{Gamma}). For $\L_\nu \eta\in\Gamma^\nu_j\cap\left\{2^{j-1}\leq|\eta|\leq2^{j+1}\right\}$, we have 
\bel{eta norm}
2^{j-1}~\leq~|\eta_\imath|~\leq~2^{j+1},\qquad |\eta_\imath^\dagger|~\leq~\C~ 2^{j/2}.
\eeq
Let $\Lambda_\ell$ be defined in (\ref{Cone}). For $\xi=(\tau,\lambda)\in \Lambda_\ell \cap\left\{2^{j-1}\leq|\xi|\leq2^{j+1}\right\}$, we have 
\bel{tau lambda norm}
\begin{array}{lr}\ds
2^{j-1}~\leq~|\lambda|~\leq~2^{j+1},\qquad |\tau|~\leq~\C ~2^{j-\ell},\qquad \ell\ge0,
\\\\ \ds
2^{j-1}~\leq~|\tau|~\leq~2^{j+1},\qquad |\lambda|~\leq~\C ~2^{j+\ell},\qquad \ell\leq0.
\end{array}
\eeq
Essentially, we require
\bel{Intersection nonempty}
\Gamma^\nu_j~\cap~\Lambda_\ell~\cap~\left\{2^{j-1}\leq|\xi|=|\eta|\leq2^{j+1}\right\}~\neq~\emptyset.
\eeq
Let $\Psi(x,\eta)$ be defined in (\ref{Phi split}). Recall from {\bf 4.5}, chapter IX of Stein \cite{Stein}. We have
\bel{d Est Psi}
\left|\left({\p\over \p\eta_\imath}\right)^\alpha\Psi(x,\eta)\right|~\leq~\C_\alpha~2^{-\alpha j},\qquad
\left|\left({\p\over\p\eta_\imath^\dagger}\right)^\beta\Psi(x,\eta)\right|~\leq~\C_\beta~2^{-|\beta| \left({1\over 2}\right)j}
\eeq
for every multi-indices $\alpha, \beta$ whenever $2^{j-1}\leq|\eta|\leq2^{j+1}$.

Let $\delta_\ell(\xi)=\delta_\ell(\L_\nu\eta)$ be defined in (\ref{delta_t}) and $\sigma(x,y,\xi)=\sigma(x,y,\L_\nu\eta)\in\S^{-{n-1\over 2}}$ satisfying the  differential inequality in (\ref{Class}).  
\begin{remark} Every $\p_\tau$ acting on $\delta_\ell(\xi)\sigma(x,y,\xi)$ gains a factor of $2^\ell|\xi|^{-1}=\C2^\ell|\lambda|^{-1}=\C|\tau|^{-1}$ for $\ell\ge0$ and every $\p_\lambda$ acting on $\delta_\ell(\xi)\sigma(x,y,\xi)$ gains a factor of $2^{-\ell}|\xi|^{-1}=\C2^{-\ell}|\tau|^{-1}=\C|\lambda|^{-1}$ for $\ell\leq0$. 
\end{remark}
By applying chain rule of differentiation, for every multi-indices $\alpha, \beta$, we have
\bel{Diff Ineq sigma eta_imath^dagger}
\begin{array}{lr}\ds
\left|\left({\p\over \p \eta_\imath^\dagger}\right)^\alpha \delta_\ell(\L_\nu\eta)\sigma(x,y,\L_\nu\eta)\right|~\leq~\C_\alpha~
\left({1\over 1+|\eta|}\right)^{n-1\over 2}2^{|\alpha|\ell}\left({1\over 1+|\eta|}\right)^{|\alpha|},\qquad\ell\ge0,
\\\\ \ds
\left|\left({\p\over \p \eta_\imath^\dagger}\right)^\beta \delta_\ell(\L_\nu\eta)\sigma(x,y,\L_\nu\eta)\right|~\leq~\C_\beta~
\left({1\over 1+|\eta|}\right)^{n-1\over 2}2^{-|\beta|\ell}\left({1\over 1+|\eta|}\right)^{|\beta|},\qquad\ell\leq0.
\end{array}
\eeq

{\bf Case One:} Let $-{j/2}-3\leq\ell\leq {j/ 2}+3$.  Note that $\xi=\L_\nu \eta$ where $\det\L_\nu=1$. We have
\bel{tau lambda entry}
\begin{array}{cc}\ds
\tau~=~a_{ 1 \imath} \eta_\imath~+~\O(1)\cdot \eta_\imath^\dagger,\qquad ( \xi_1=\tau )
\\\\ \ds
\lambda_i~=~a_{i\imath}\eta_\imath~+~\O(1)\cdot\eta_\imath^\dagger, \qquad i=2,\ldots,n
\end{array}
\eeq
where  $a_{i\imath}$ denote the entry on the  $i$-th row and $\imath$-th column of $\L_\nu$.  

By putting together (\ref{eta norm})-(\ref{tau lambda norm}) and (\ref{tau lambda entry}), we necessarily  have
\bel{entry norm}
\begin{array}{cc}\ds
|a_{1\imath}|~\leq~\C2^{-\ell}, \qquad 0~\leq~\ell~\leq~ j/2+3,
\\\\ \ds
|a_{i\imath}|\leq\C2^{\ell}, \qquad i~=~2,\ldots,n, \qquad -j/2-3~\leq~\ell~\leq~0.
\end{array}
\eeq
Recall from {\bf Remark 7.1}. By  applying chain rule of differentiation and using  (\ref{entry norm}), we have
\bel{Diff Ineq sigma eta_imath}
\begin{array}{rl}\ds
\left|\left({\p\over \p \eta_\imath}\right)^\alpha \delta_\ell(\L_\nu\eta)\sigma(x,y,\L_\nu\eta)\right|~\leq~\C_\alpha~\left({1\over 1+|\eta|}\right)^{n-1\over 2}\left({1\over 1+|\eta|}\right)^\alpha
\end{array}
\eeq
for every  $\alpha$.

Define 
\bel{L operator}
\mathcal{L}~=~I+2^{2j}\left({\p\over \p\eta_\imath}\right)^2+2^j\Delta_{\eta_\imath^\dagger}
\eeq
where $I$ is the identity operator. 

Let $\Theta^\nu_{\ell j}(x,y,\eta)$ be defined in (\ref{Omega rewrite}).
Recall the estimate in (\ref{chi prop est split}), together with (\ref{Theta norm}), (\ref{d Est Psi})  and (\ref{Diff Ineq sigma eta_imath^dagger})-(\ref{Diff Ineq sigma eta_imath}).  Note that $\ell\leq j/2+3$. We have
\bel{Theta diff est}
\left|{\mathcal{L}}^N \Theta^\nu_{\ell j}(x,y,\eta)\right|~\leq~\C_N~2^{-j\left({n-1\over 2}\right)},\qquad N\ge1.
\eeq
Moreover, by (\ref{intersection of cone norm}), the support of $\Theta^\nu_{\ell j}(x,y,\eta)$ in $\eta$ has a volume bounded by $\C~2^j2^{j\left({n-1\over 2}\right)}$.

Let $\Omega^\nu_{\ell j}(x,y)$ be given in (\ref{Omega rewrite}).  By using (\ref{Theta diff est}), an $N$-fold integration by parts associated to $\mathcal{L}$ shows that
\bel{Omega rewrite norm}
\begin{array}{lr}\ds
\left|\Omega^\nu_{\ell j}(x,y)\right|~\leq~\C_{N}~2^{-j\left({n-1\over 2}\right)} 2^j2^{j\left({n-1\over 2}\right)}
\\\\ \ds
\left\{1+4\pi^2 2^{2j}\left(\nabla_\eta\Phi\left(x,\L_\nu\eta_j^\nu\right)-\L_\nu^T y\right)_\imath^2+4\pi^2 2^{j}    
\left| \left(\nabla_\eta\Phi\left(x,\L_\nu\eta_j^\nu\right)-\L_\nu^T y\right)_\imath^\dagger \right|^2\right\}^{-N}. 
\end{array}
\eeq

Consider  the local diffeomorphism  $\mathcal{X}_\Phi~\colon x\mt\left(\L_\nu^T\right)^{-1}\nabla_\eta\Phi\left(x,\L_\nu\eta_j^\nu\right)$
whose Jacobian has an absolute value strictly greater than zero. Indeed, recall from (\ref{det Phi_xxi}) where we have $\det\Phi_{x\xi}(x,\xi)=\det\Phi_{x\eta}(x,\L_\nu \eta)\neq0$ for $\eta\ne0$ provided that  the nondegeneracy condition hold in (\ref{nondegeneracy}). 

Denote $\mathcal{X}=\mathcal{X}(x)=\left(\L_\nu^T\right)^{-1}\nabla_\eta\Phi\left(x,\L_\nu\eta_j^\nu\right)$. By using  (\ref{Omega rewrite norm}), we have
\bel{int Omega rewrite norm}
\begin{array}{lr}\ds
\int_{\R^n}\left|\Omega^\nu_{\ell j}(x,y)\right|dx
\\\\ \ds
~\leq~\C_{\Phi~N}\int_{\R^n} 2^{-j\left({n-1\over 2}\right)}2^j2^{j\left({n-1\over 2}\right)}\left\{1+ 2^{2j}(\mathcal{X}-y)_\imath^2+ 2^{j}    
\left| (\mathcal{X}-y)_\imath^\dagger \right|^2\right\}^{-N} d\mathcal{X}
\\\\ \ds
~\leq~\C_{\Phi~N}\int_{\R^n} 2^{-j\left({n-1\over 2}\right)}\left\{1+\mathcal{Z}_\imath^2+     
\left| \mathcal{Z}_\imath^\dagger\right|^2\right\}^{-N} d\mathcal{Z} 
\\\\ \ds~~~~~~~~~~~~~~~~~~~~~~~~~~~~~~~~~~~~~~~~~~~~~~~~~~~
\hbox{\small{$ \mathcal{Z}_\imath=2^{j}(\mathcal{X}-y)_\imath,~~~~ \mathcal{Z}_\imath^\dagger=2^{j/2}(\mathcal{X}-y)_\imath^\dagger$}}
\\\\ \ds
~\leq~\C_{\Phi}~2^{-j\left({n-1\over 2}\right)},\qquad \hbox{\small{$N\ge {n+1\over 2}$}}.
\end{array}
\eeq
Recall from (\ref{number elements case3}). For $0\leq\ell\leq j/2+3$,  the dyadic cone $\Lambda_\ell$ intersects at most  $\C ~2^{j\left({n-1\over 2}\right)}2^{-\ell}$ many of $\Gamma^\nu_j$. On the other hand, for $-j/2-3\leq\ell\leq0$,  the dyadic cone $\Lambda_\ell$ intersects at most  $\C ~2^{j\left({n-1\over 2}\right)}2^{\ell(n-1)}$ many of $\Gamma^\nu_j$ .

We simultaneously have
\bel{Omega_lj Sum norm >}
\begin{array}{lr}\ds
\int_{\R^n}\left|\Omega_{\ell j}(x,y)\right|dx~\leq~\sum_\nu \int_{\R^n}\left|\Omega^\nu_{\ell j}(x,y)\right|dx\qquad\hbox{\small{by (\ref{Omega_lj Sum})}}
\\\\ \ds ~~~~~~~~~~~~~~~~~~~~~~~~~~~~
~\leq~\C_{\Phi}~ 2^{j\left({n-1\over 2}\right)}2^{-\ell}~2^{-j\left({n-1\over 2}\right)}\qquad \hbox{\small{by (\ref{int Omega rewrite norm})}}
\\\\ \ds ~~~~~~~~~~~~~~~~~~~~~~~~~~~~
~\leq~\C_{\Phi}~2^{-\ell},\qquad \hbox{\small{$0\leq\ell\leq j/2+3$}}
\end{array}
\eeq
and
\bel{Omega_lj Sum norm <}
\begin{array}{lr}\ds
\int_{\R^n}\left|\Omega_{\ell j}(x,y)\right|dx~\leq~\sum_\nu \int_{\R^n}\left|\Omega^\nu_{\ell j}(x,y)\right|dx
\qquad\hbox{\small{by (\ref{Omega_lj Sum})}}
\\\\ \ds ~~~~~~~~~~~~~~~~~~~~~~~~~~~~
~\leq~\C_{\Phi}~ 2^{j\left({n-1\over 2}\right)}2^{\ell(n-1)}~2^{-j\left({n-1\over 2}\right)}\qquad \hbox{\small{by (\ref{int Omega rewrite norm})}}
\\\\ \ds ~~~~~~~~~~~~~~~~~~~~~~~~~~~~
~\leq~\C_{\Phi}~2^{\ell(n-1)},\qquad \hbox{\small{$- j/2-3\leq\ell\leq0$.}}
\end{array}
\eeq
Observe that every $\p_y$ acting on $\Omega^\nu_{\ell j}(x,y)$ defined in (\ref{Omega^v_lj}) or (\ref{Omega rewrite}) gains a factor of $\C 2^j$ whenever $2^{j-1}\leq|\xi|=|\eta|\leq2^{j+1}$. By carrying out the same estimates in (\ref{Theta norm})-(\ref{Omega_lj Sum norm <}), we find
\bel{nabla Omega_lj} 
\int_{\R^n} \left|\nabla_y\Omega_{\ell j}(x,y)\right| dx~\leq~\C_{\Phi}~
 \left\{\begin{array}{lr}\ds 2^j2^{-\ell},\qquad ~~~~~0\leq\ell\leq j/2+3,
 \\\\ \ds
2^j 2^{(n-1)\ell},\qquad -j/2-3\leq\ell\leq0
 \end{array}\right. 
 \eeq
which further implies 
\bel{Est2 =} 
\int_{\R^n} \left|\Omega_{\ell j}(x,y)-\Omega_{\ell j}(x,x_o)\right| dx~\leq~\C_{\Phi}~2^j|y-x_o|
 \left\{\begin{array}{lr}\ds 2^{-\ell},\qquad ~~~~~0\leq\ell\leq j/2+3,
 \\\\ \ds
 2^{(n-1)\ell},\qquad -j/2-3\leq\ell\leq0.
 \end{array}\right. 
 \eeq

Recall $\mathfrak{Q}_r(x_o)$ defined in (\ref{rectangle R})-(\ref{Q_r}). Let $2^{k}\leq r^{-1}\leq 2^{k+1}$.
For $x\in{^c}\mathfrak{Q}_r(x_o)$, we must  have
\bel{x-Phi Est =}
\left|\left(\L_\nu^Tx_o-\nabla_\eta\Phi\left(x,\L_\nu \eta^\nu_j\right)\right)_\imath\right|~\ge~2\cdot2^{-k}~~~~\hbox{or}~~~~
\left|\left(\L_\nu^Tx_o-\nabla_\eta\Phi\left(x,\L_\nu\eta^\nu_j\right)\right)_\imath^\dagger\right|~\ge~ 2\cdot2^{-k/2}.
\eeq
If $y\in B_r(x_o)$, then $|y-x_o|\leq 2^{-k}$. For every $2^j\ge r^{-1}$,  we have
\bel{x,y-Phi Est =}
2^{2j}\left|\left(\L_\nu^Ty-\nabla_\eta\Phi\left(x,\L_\nu \eta^\nu_j\right)\right)_\imath\right|^2+
2^j\left|\left(\L_\nu^Ty-\nabla_\eta\Phi\left(x,\L_\nu\eta^\nu_j\right)\right)_\imath^\dagger\right|^2~\ge~ 2^{j-k}.
\eeq
By carrying out the same estimates in (\ref{Theta norm})-(\ref{Omega_lj Sum norm <}), except for (\ref{Omega rewrite norm})  replaced with 
\bel{Omega rewrite norm*}
\begin{array}{lr}\ds
\left|\Omega^\nu_{\ell j}(x,y)\right|~\leq~\C_{N}~2^{-j\left({n-1\over 2}\right)} 2^j2^{j\left({n-1\over 2}\right)}2^{-j+k}
\\\\ \ds
\left\{1+4\pi^2 2^{2j}\left(\nabla_\eta\Phi\left(x,\L_\nu\eta_j^\nu\right)-\L_\nu^T y\right)_\imath^2+4\pi^2 2^{j}    
\left| \left(\nabla_\eta\Phi\left(x,\L_\nu\eta_j^\nu\right)-\L_\nu^T y\right)_\imath^\dagger \right|^2\right\}^{1-N} 
~~~\hbox{\small{by (\ref{x,y-Phi Est =})}}

\end{array}
\eeq
 we find
\bel{Est3 =} 
\int_{{^c}\mathfrak{Q}_r(x_o)} \left|\Omega_{\ell j}(x,y)\right| dx~\leq~\C_{\Phi}~ {2^{-j}\over r}
 \left\{\begin{array}{lr}\ds 2^{-\ell},\qquad ~~~~~0\leq\ell\leq j/2+3,
 \\\\ \ds
 2^{(n-1)\ell},\qquad -j/2-3\leq\ell\leq0
 \end{array}\right. 
 \eeq
for every $y\in B_r(x_o)$ whenever $2^{j}\ge r^{-1}$.

\v

{\bf Case Two:} Let $\ell>j/2+3$. Recall from section 6. We consider $\xi^\nu_j$ equally distributed on the unit sphere $\mathbb{S}^{n-2}$ of the $(n-1)$-dimensional $\lambda$-space.   $\L_\nu$ defined in (\ref{L_v decomposable}) is decomposable:
\bel{L_v decomposable ^}
\begin{array}{cc}\ds
\L_\nu~=~\left[\begin{array}{lr}
1 ~~~~
\\ \ds
~~~~ \L'_\nu
\end{array}\right],\qquad \det\L'_\nu~=~1
\\\\ \ds
\hbox{and}\qquad \xi^\nu_j~=~\L_\nu \eta^\nu_j,\qquad \eta^\nu_j~=~\left({\eta_\imath\over|\eta_\imath|},0\right)~=~\left({\eta_n\over|\eta_n|},0\right).
 \end{array}
\eeq
In particular,   $\xi_1=\tau=\eta_1$ is independent from $\eta_\imath=\eta_n$.

Let $\delta_\ell(\xi)=\delta(\L_\nu\eta)$ be defined in (\ref{delta_t}) and $\sigma(x,y,\xi)=\sigma(x,y,\L_\nu\eta)\in\S^{-{n-1\over 2}}$ satisfying the differential inequality in (\ref{Class}). We have
\bel{Diff Ineq sigma eta_imath >}
\begin{array}{lr}\ds
\left|\left({\p\over \p \eta_n}\right)^\alpha \delta_\ell(\L_\nu\eta)\sigma(x,y,\L_\nu\eta)\right|~\leq~\C_\alpha~\left({1\over 1+|\eta|}\right)^{n-1\over 2}\left({1\over 1+|\eta|}\right)^\alpha
\end{array}
\eeq
for every  $\alpha$.

Suppose $\ell>j+3$.
Define 
\bel{D operator sharp}
\mathcal{D^\sharp}~=~I+2^{2j}\left({\p\over \p\eta_n}\right)^2+2^j\sum_{i=2}^{n-1}\left({\p\over \p\eta_i}\right)^2.
\eeq
Let $\Theta^\nu_{\ell j}(x,y,\eta)$ be defined in (\ref{Omega rewrite}).
Recall the estimate in (\ref{chi prop est split}), together with (\ref{Theta norm}), (\ref{d Est Psi})  and (\ref{Diff Ineq sigma eta_imath >}).  We have
\bel{Theta diff est > sharp}
\left|\left(\mathcal{D}^\sharp\right)^N \Theta^\nu_{\ell j}(x,y,\eta)\right|~\leq~\C_N~2^{-j\left({n-1\over 2}\right)},\qquad N\ge1.
\eeq
Moreover, by (\ref{intersection of cone norm>}), the support of $\Theta^\nu_{\ell j}(x,y,\eta)$ in $\eta$ has a volume bounded by $\C~2^{j-\ell}2^j2^{j\left({n-2\over 2}\right)}$.

Let $\Omega^\nu_{\ell j}(x,y)$ be defined in (\ref{Omega rewrite}). By using (\ref{Theta diff est > sharp}), an $N$-fold integration by parts associated to $\mathcal{D}^\sharp$ shows that
\bel{Omega rewrite norm > sharp}
\begin{array}{lr}\ds
\left|\Omega^\nu_{\ell j}(x,y)\right|~\leq~\C_{N}~2^{-j\left({n-1\over 2}\right)}~2^{j-\ell}2^j2^{j\left({n-2\over 2}\right)}
\\\\ \ds
\left\{1+4\pi^2 2^{2j}\left(\nabla_\eta\Phi\left(x,\L_\nu\eta_j^\nu\right)-\L_\nu^T y\right)_n^2+4\pi^2 2^{j}\sum_{i=2}^{n-1}\left(\nabla_\eta\Phi\left(x,\L_\nu\eta_j^\nu\right)-\L_\nu^T y\right)_i^2 \right\}^{-N} . 
\end{array}
\eeq
Note that $\Omega^\nu_{\ell j}(x,y)$ has the same compact support of $\sigma(x,y,\xi)$ in both $x$ and $y$.

Recall  the local diffeomorphism  $\mathcal{X}_\Phi\colon x\mt\left(\L_\nu^T\right)^{-1}\nabla_\eta\Phi\left(x,\L_\nu\eta_j^\nu\right)$ whose Jacobian has an absolute value strictly greater than zero. Let $\mathcal{X}_\Phi \supp\sigma$ denote the image of the $x$-support of $\sigma(x,y,\xi)$  and
$\mathcal{X}=\mathcal{X}(x)=\left(\L_\nu^T\right)^{-1}\nabla_\eta\Phi\left(x,\L_\nu\eta_j^\nu\right)$.
By  using (\ref{Omega rewrite norm > sharp}), we have
\bel{int Omega rewrite norm > sharp}
\begin{array}{lr}\ds
\int_{\R^n}\left|\Omega^\nu_{\ell j}(x,y)\right|dx
\\\\ \ds
~\leq~\C_{\Phi~N}\int_{\mathcal{X}_\Phi\supp\sigma} 2^{-j\left({n-1\over 2}\right)}2^{j-\ell}2^j2^{j\left({n-2\over 2}\right)}\left\{1+ 2^{2j}(\mathcal{X}-y)_n^2+2^j\sum_{i=2}^{n-1}(\mathcal{X}-y)_i^2\right\}^{-N}d\mathcal{X} 
\\\\ \ds
~\leq~\C_{\Phi~N}\iint_{\R^{n-1} ~\times~ \mathcal{X}_\Phi\supp\sigma\cap\R} 2^{-j\left({n-1\over 2}\right)}2^{j-\ell}\left\{1+ \mathcal{Z}_n^2+     \sum_{i=2}^{n-1} \mathcal{Z}_i^2 \right\}^{-N} d\mathcal{Z}_1^\dagger d\mathcal{Z}_1 
\\\\ \ds~~~~~~~~~~~~~~~~~~~~~~~~~~~~~~~~~~~~~~
\hbox{\small{$ \mathcal{Z}_n=2^{j}(\mathcal{X}-y)_n$~~and~~ $\mathcal{Z}_i=2^{j/2}(\mathcal{X}-y)_i,~~i=2,\ldots,n-1$}}
\\\\ \ds
~\leq~\C_{\sigma~\Phi}~2^{-j\left({n-1\over 2}\right)}2^{j-\ell},\qquad N\ge {n/2}.
\end{array}
\eeq

Recall from section 6. For $\ell>j/2+3$, $\Lambda_\ell$ intersects  every $\Gamma^\nu_j$ whose central direction is $\xi^\nu_j\in\mathbb{S}^{n-2}$ in the $\lambda$-space.  There are at most $\C ~2^{\left({n-2\over 2}\right)j}$ such elements in $\left\{\xi^\nu_j\right\}_\nu$.

We thus have
\bel{Omega_lj Sum norm >sharp}
\begin{array}{lr}\ds
\int_{\R^n}\left|\Omega_{\ell j}(x,y)\right|dx~\leq~\sum_{\nu~\colon~\xi^\nu_j\in\mathbb{S}^{n-2}} \int_{\R^n}\left|\Omega^\nu_{\ell j}(x,y)\right|dx\qquad \hbox{\small{by (\ref{Omega_lj Sum})}}
\\\\ \ds ~~~~~~~~~~~~~~~~~~~~~~~~~~~~
~\leq~\C_{\sigma~\Phi}~ 2^{j\left({n-2\over 2}\right)} 2^{-j\left({n-1\over 2}\right)}2^{j-\ell}~=~\C_{\sigma~\Phi}~2^{\left({1\over 2}\right)j-\ell}
\qquad \hbox{\small{by (\ref{int Omega rewrite norm > sharp})}}
\\\\ \ds ~~~~~~~~~~~~~~~~~~~~~~~~~~~~
~\leq~\C_{\sigma~\Phi}~2^{\left({1\over 2}\right)\ell-\ell},\qquad (~\ell>j+3~)
\\\\ \ds ~~~~~~~~~~~~~~~~~~~~~~~~~~~~
~\leq~\C_{\sigma~\Phi}~2^{-\left({1\over 2}\right)\ell}.
\end{array}
\eeq

Suppose $j/2+3<\ell\leq j+3$.  Define 
\bel{D operator flat}
\mathcal{D^\flat}~=~I+2^{2j}\left({\p\over \p\eta_n}\right)^2+2^j\sum_{i=2}^{n-1}\left({\p\over \p\eta_i}\right)^2+2^{2(j-\ell)}\left({\p\over \p\eta_1}\right)^2.
\eeq
Let $\Theta^\nu_{\ell j}(x,y,\eta)$ be defined in (\ref{Omega rewrite}).
Recall the estimate in (\ref{chi prop est split}), together with (\ref{Theta norm}), (\ref{d Est Psi})  and (\ref{Diff Ineq sigma eta_imath >}).  Note that  $\eta_1=\tau$ and every $\p_\tau$ acting on  
$\chi^\nu_{\ell j}(\L_\nu\eta)\delta_\ell(\L_\nu\eta)\phi_j(\L_\nu\eta)\sigma(x,y,\L_\nu\eta)$ gains a factor bounded by $\C2^{-j+\ell}$.
We have
\bel{Theta diff est > flat}
\left|\left(\mathcal{D}^\flat\right)^N \Theta^\nu_{\ell j}(x,y,\eta)\right|~\leq~\C_N~2^{-j\left({n-1\over 2}\right)},\qquad N\ge1.
\eeq
Moreover,  by (\ref{intersection of cone norm>}), the support of $\Theta^\nu_{\ell j}(x,y,\eta)$ in $\eta$ has a volume bounded by $\C~2^{j-\ell}2^j2^{j\left({n-2\over 2}\right)}$.

Let $\Omega^\nu_{\ell j}(x,y)$ be given in (\ref{Omega rewrite}). By using (\ref{Theta diff est > flat}), an $N$-fold integration by parts associated to $\mathcal{D}^\flat$ shows that
\bel{Omega rewrite norm > flat}
\begin{array}{lr}\ds
\left|\Omega^\nu_{\ell j}(x,y)\right|~\leq~\C_{N} ~2^{-j\left({n-1\over 2}\right)}~2^{j-\ell}2^j2^{j\left({n-2\over 2}\right)}
\\\\ \ds
\Bigg\{1+4\pi^2 2^{2j}\left(\nabla_\eta\Phi\left(x,\L_\nu\eta_j^\nu\right)-\L_\nu^T y\right)_n^2+4\pi^2 2^{j}\sum_{i=2}^{n-1}\left(\nabla_\eta\Phi\left(x,\L_\nu\eta_j^\nu\right)-\L_\nu^T y\right)_i^2 
\\\\ \ds
+4\pi^2 2^{2(j-\ell)}\left(\nabla_\eta\Phi\left(x,\L_\nu\eta_j^\nu\right)-\L_\nu^T y\right)_1^2\Bigg\}^{-N}. 
\end{array}
\eeq
Denote $\mathcal{X}=\mathcal{X}(x)=\left(\L_\nu^T\right)^{-1}\nabla_\eta\Phi\left(x,\L_\nu\eta_j^\nu\right)$.
By using (\ref{Omega rewrite norm > flat}), we have
\bel{int Omega rewrite norm > flat}
\begin{array}{lr}\ds
\int_{\R^n}\left|\Omega^\nu_{\ell j}(x,y)\right|dx
\\\\ \ds
~\leq~\C_{\Phi~N}\int_{\R^n} 2^{-j\left({n-1\over 2}\right)}2^{j-\ell}2^j2^{j\left({n-2\over 2}\right)}\left\{1+ 2^{2j}(\mathcal{X}-y)_n^2+2^j\sum_{i=2}^{n-1}(\mathcal{X}-y)_i^2+2^{2(j-\ell)}(\mathcal{X}-y)_1^2\right\}^{-N}d\mathcal{X} 
\\\\ \ds
~\leq~\C_{\Phi~N}\int_{\R^n} 2^{-j\left({n-1\over 2}\right)}\left\{1+ \mathcal{Z}_n^2+     \sum_{i=2}^{n-1} \mathcal{Z}_i^2 +\mathcal{Z}_1^2\right\}^{-N} dx
\\\\ \ds~~~~~~~~~~~~~~~~~~~~~~~~~
\hbox{\small{$ \mathcal{Z}_n=2^{j}(\mathcal{X}-y)_n,~~\mathcal{Z}_1=2^{j-\ell}(\mathcal{X}-y)_1$}} 
~~
\hbox{\small{ and ~~$ \mathcal{Z}_i=2^{j/2}(\mathcal{X}-y)_i,~~i=2,\ldots,n-1$}}
\\\\ \ds
~\leq~\C_{\sigma~\Phi}~2^{-j\left({n-1\over 2}\right)},\qquad N\ge {(n+1)/2}.
\end{array}
\eeq
For $\ell>j/2+3$, $\Lambda_\ell$ intersects  every $\Gamma^\nu_j$ whose central direction is $\xi^\nu_j\in\mathbb{S}^{n-2}$ in the $\lambda$-space.  There are at most $\C ~2^{\left({n-2\over 2}\right)j}$ such elements in $\left\{\xi^{\nu}_j\right\}_\nu$.
We thus have
\bel{Omega_lj Sum norm>flat}
\begin{array}{lr}\ds
\int_{\R^n}\left|\Omega_{\ell j}(x,y)\right|dx~\leq~\sum_{\nu~\colon~\xi^\nu_j\in\mathbb{S}^{n-2}} \int_{\R^n}\left|\Omega^\nu_{\ell j}(x,y)\right|dx\qquad \hbox{\small{by (\ref{Omega_lj Sum})}}
\\\\ \ds ~~~~~~~~~~~~~~~~~~~~~~~~~~~~
~\leq~\C_{\sigma~\Phi}~ 2^{j\left({n-2\over 2}\right)} 2^{-j\left({n-1\over 2}\right)}~=~\C_{\sigma~\Phi}~2^{-\left({1\over 2}\right)j}
\qquad \hbox{\small{by (\ref{int Omega rewrite norm > flat})}}
\\\\ \ds ~~~~~~~~~~~~~~~~~~~~~~~~~~~~
~\leq~\C_{\sigma~\Phi}~2^{-\left({1\over 2}\right)\ell}.\qquad (~\ell\leq j+3~)
\end{array}
\eeq
Observe that every $\p_y$ acting on $\Omega^\nu_{\ell j}(x,y)$ defined in (\ref{Omega^v_lj}) or (\ref{Omega rewrite}) gains a factor of $\C 2^j$ whenever $2^{j-1}\leq|\xi|=|\eta|\leq2^{j+1}$. By carrying out the same estimates in (\ref{L_v decomposable ^})-(\ref{Omega_lj Sum norm>flat}), we find
\bel{nabla Omega_lj>} 
\int_{\R^n} \left|\nabla_y\Omega_{\ell j}(x,y)\right| dx~\leq~\C_{\sigma~\Phi}~
 2^j2^{-\left({1\over 2}\right)\ell}
\eeq
which further implies 
\bel{Est2 >} 
\int_{\R^n} \left|\Omega_{\ell j}(x,y)-\Omega_{\ell j}(x,x_o)\right| dx~\leq~\C_{\sigma~\Phi}~2^j|y-x_o|
 2^{-\left({1\over 2}\right)\ell}.
  \eeq
\v

Recall $\mathfrak{Q}_r(x_o)$ defined in (\ref{rectangle R>})-(\ref{Q_r>}). Let $2^{k}\leq r^{-1}\leq 2^{k+1}$.
For $x\in{^c}\mathfrak{Q}_r(x_o)$, we must have
\bel{x-Phi Est >}
\left|\left(\L_\nu^Tx_o-\nabla_\eta\Phi\left(x,\L_\nu \eta^\nu_j\right)\right)_n\right|~\ge~2\cdot2^{-k}~~~~\hbox{or}~~~~
\left|\left(\L_\nu^Tx_o-\nabla_\eta\Phi\left(x,\L_\nu\eta^\nu_j\right)\right)_1^\dagger\right|~\ge~ 2\cdot2^{-k/2}.
\eeq
If $y\in B_r(x_o)$, then $|y-x_o|\leq 2^{-k}$. For every $2^j\ge r^{-1}$,  we have
\bel{x,y-Phi Est >}
2^{2j}\left(\L_\nu^Ty-\nabla_\eta\Phi\left(x,\L_\nu \eta^\nu_j\right)\right)_n^2+
2^j\sum_{i=2}^n \left(\L_\nu^Ty-\nabla_\eta\Phi\left(x,\L_\nu\eta^\nu_j\right)\right)_i^2~\ge~ 2^{j-k}.
\eeq
By carrying out the same estimates in (\ref{L_v decomposable ^})-(\ref{Omega_lj Sum norm>flat}), except for (\ref{Omega rewrite norm > sharp}) and (\ref{Omega rewrite norm > flat}) replaced  with
\bel{Omega rewrite norm > sharp*}
\begin{array}{lr}\ds
\left|\Omega^\nu_{\ell j}(x,y)\right|~\leq~\C_{N}~2^{-j\left({n-1\over 2}\right)}~2^{j-\ell}2^j2^{j\left({n-2\over 2}\right)}
~2^{-j+k}
\\\\ \ds
\left\{1+4\pi^2 2^{2j}\left(\nabla_\eta\Phi\left(x,\L_\nu\eta_j^\nu\right)-\L_\nu^T y\right)_n^2+4\pi^2 2^{j}\sum_{i=2}^{n-1}\left(\nabla_\eta\Phi\left(x,\L_\nu\eta_j^\nu\right)-\L_\nu^T y\right)_i^2 \right\}^{1-N} 
\qquad\hbox{\small{by (\ref{x,y-Phi Est >})}}
\end{array}
\eeq
 and
\bel{Omega rewrite norm > flat*}
\begin{array}{lr}\ds
\left|\Omega^\nu_{\ell j}(x,y)\right|~\leq~\C_{N} ~2^{-j\left({n-1\over 2}\right)}~2^{j-\ell}2^j2^{j\left({n-2\over 2}\right)}
~2^{-j+k}
\\\\ \ds
\Bigg\{1+4\pi^2 2^{2j}\left(\nabla_\eta\Phi\left(x,\L_\nu\eta_j^\nu\right)-\L_\nu^T y\right)_n^2+4\pi^2 2^{j}\sum_{i=2}^{n-2}\left(\nabla_\eta\Phi\left(x,\L_\nu\eta_j^\nu\right)-\L_\nu^T y\right)_i^2 
\\\\ \ds 
+4\pi^2 2^{2(j-\ell)}\left(\nabla_\eta\Phi\left(x,\L_\nu\eta_j^\nu\right)-\L_\nu^T y\right)_1^2\Bigg\}^{1-N}~~~~~~~~\hbox{\small{by (\ref{x,y-Phi Est >})}}
\end{array}
\eeq
 we find
\bel{Est3 >} 
\int_{{^c}\mathfrak{Q}_r(x_o)} \left|\Omega_{\ell j}(x,y)\right| dx~\leq~\C_{\sigma~\Phi}~ {2^{-j}\over r}~
  2^{-\left({1\over 2}\right)\ell} 
 \eeq
for every $y\in B_r(x_o)$ whenever $2^{j}\ge r^{-1}$.

\v

{\bf Case Three:} Let $\ell<-j/2-3$. Recall from section 6. We consider  only for $\xi^\nu_j=(\tau,\lambda)=(\pm1,0)$.
$\L_\nu$ is  the identity matrix 
so that $\xi=\eta$. In particular, $\eta_\imath=\eta_1=\xi_1=\tau$.

Let $\delta_\ell(\xi)=\delta_\ell(\L_\nu \eta)$ be defined in (\ref{delta_t}) and $\sigma(x,y,\xi)=\sigma(x,y,\L_\nu \eta)\in\S^{-{n-1\over 2}}$ satisfying the differential inequality in (\ref{Class}). 

We have
\bel{Diff Ineq sigma eta_imath <}
\begin{array}{lr}\ds
\left|\left({\p\over \p \eta_1}\right)^\alpha \delta_\ell(\L_\nu\eta)\sigma(x,y,\L_\nu\eta)\right|~\leq~\C_\alpha~\left({1\over 1+|\eta|}\right)^{n-1\over 2}\left({1\over 1+|\eta|}\right)^\alpha
\end{array}
\eeq
for every $\alpha$.

Suppose $\ell<-j-3$.
Define 
\bel{G operator}
\mathfrak{D^\sharp}~=~I+2^{2j}\left({\p\over \p\eta_1}\right)^2.
\eeq
Let $\Theta^\nu_{\ell j}(x,y,\eta)$ be defined in (\ref{Omega rewrite}).
Recall the estimate in (\ref{chi prop est split}), together with (\ref{Theta norm}), (\ref{d Est Psi})  and (\ref{Diff Ineq sigma eta_imath <}).  We have
\bel{Theta diff est < sharp}
\left|\left(\mathfrak{D}^\sharp\right)^N \Theta^\nu_{\ell j}(x,y,\eta)\right|~\leq~\C_N~2^{-j\left({n-1\over 2}\right)},\qquad N\ge1.
\eeq
Moreover,  by (\ref{intersection of cone norm<}), the support of $\Theta^\nu_{\ell j}(x,y,\eta)$ in $\eta$ has a volume bounded by $\C~2^j2^{(j+\ell)\left({n-1}\right)}$.

Let $\Omega^\nu_{\ell j}(x,y)$ be given in (\ref{Omega rewrite}). By using (\ref{Theta diff est < sharp}), an $N$-fold integration by parts associated to $\mathfrak{D}^\sharp$ shows that
\bel{Omega rewrite norm < sharp}
\begin{array}{lr}\ds
\left|\Omega^\nu_{\ell j}(x,y)\right|~\leq~\C_{N} ~2^{-j\left({n-1\over 2}\right)}2^j2^{(j+\ell)\left({n-1}\right)}
\left\{1+4\pi^2 2^{2j}\left(\nabla_\eta\Phi\left(x,\L_\nu\eta_j^\nu\right)-\L_\nu^T y\right)_1^2 \right\}^{-N}. 
\end{array}
\eeq
Recall  the local diffeomorphism  $\mathcal{X}_\Phi\colon x\mt\left(\L_\nu^T\right)^{-1}\nabla_\eta\Phi\left(x,\L_\nu\eta_j^\nu\right)$ whose Jacobian has an absolute value strictly greater than zero. Let $\mathcal{X}_\Phi \supp\sigma$ denote the image of the $x$-support of $\sigma(x,y,\xi)$  and
$\mathcal{X}=\mathcal{X}(x)=\left(\L_\nu^T\right)^{-1}\nabla_\eta\Phi\left(x,\L_\nu\eta_j^\nu\right)$.
By using (\ref{Omega rewrite norm < sharp}), we have
\bel{int Omega rewrite norm < sharp}
\begin{array}{lr}\ds
\int_{\R^n}\left|\Omega^\nu_{\ell j}(x,y)\right|dx
\\\\ \ds
~\leq~\C_{\Phi~N}\int_{\mathcal{X}_\Phi\supp\sigma} 2^{-j\left({n-1\over 2}\right)}2^j2^{(j+\ell)\left({n-1}\right)}
\left\{1+ 2^{2j}(\mathcal{X}-y)_1^2\right\}^{-N}d\mathcal{X} 
\\\\ \ds
~\leq~\C_{\Phi~N}\iint_{\R ~\times~ \mathcal{X}_\Phi\supp\sigma\cap\R^{n-1}} 2^{-j\left({n-1\over 2}\right)}2^{(j+\ell)(n-1)}\left\{1+ \mathcal{Z}_1^2 \right\}^{-N} d\mathcal{Z}_1 d\mathcal{Z}_1^\dagger 
\\\\ \ds~~~~~~~~~~~~~~~~~~~~~~~~~~~~~~~~~~~~~~~~~~~~~~~~~~~~~~~~~~~~~~~~~
\hbox{\small{$ \mathcal{Z}_1=2^{j}(\mathcal{X}-y)_1$}}
\\\\ \ds
~\leq~\C_{\sigma~\Phi}~2^{-j\left({n-1\over 2}\right)}2^{(j+\ell)(n-1)},\qquad N\ge1.
\end{array}
\eeq
Recall from section 6. For $\ell<-j/2-3$, $\Lambda_\ell$ intersects  only $\Gamma^\nu_j$ whose central direction is  $\xi^\nu_j=(\pm1,0)\in\R\times\R^{n-1}$.
We thus have
\bel{Omega_lj Sum norm < sharp}
\begin{array}{lr}\ds
\int_{\R^n}\left|\Omega_{\ell j}(x,y)\right|dx~\leq~\sum_{\nu~\colon~\xi^\nu_j=(\pm1,0)} \int_{\R^n}\left|\Omega^\nu_{\ell j}(x,y)\right|dx\qquad\hbox{\small{by (\ref{Omega_lj Sum})}}
\\\\ \ds ~~~~~~~~~~~~~~~~~~~~~~~~~~~~
~\leq~\C_{\sigma~\Phi}~ 2^{-j\left({n-1\over 2}\right)}2^{(j+\ell)(n-1)}~=~\C_{\sigma~\Phi}~2^{\left({n-1\over 2}\right)j+\ell(n-1)}\qquad \hbox{\small{by (\ref{int Omega rewrite norm < sharp})}}
\\\\ \ds ~~~~~~~~~~~~~~~~~~~~~~~~~~~~
~\leq~\C_{\sigma~\Phi}~2^{-\left({n-1\over 2}\right)\ell+\ell(n-1)},\qquad (~\ell<-j-3~)
\\\\ \ds ~~~~~~~~~~~~~~~~~~~~~~~~~~~~
~\leq~\C_{\sigma~\Phi}~2^{\left({n-1\over 2}\right)\ell}.
\end{array}
\eeq

Suppose $ -j-3\leq\ell<-j/2-3$.
Define 
\bel{G operator flat}
\mathfrak{D^\flat}~=~I+2^{2j}\left({\p\over \p\eta_1}\right)^2+2^{2(j+\ell)}\Delta_{\eta_1^\dagger}.
\eeq
Let $\Theta^\nu_{\ell j}(x,y,\eta)$ be defined in (\ref{Omega rewrite}).
Recall the estimate in (\ref{chi prop est split}), together with (\ref{Theta norm}), (\ref{d Est Psi}) and (\ref{Diff Ineq sigma eta_imath <}).  Note that  $\eta_1^\dagger=\lambda$ and every $\p_\lambda$ acting on  
$\chi^\nu_{\ell j}(\L_\nu\eta)\delta_\ell(\L_\nu\eta)\phi_j(\L_\nu\eta)\sigma(x,y,\L_\nu\eta)$ gains a factor bounded by $\C2^{-j-\ell}$.

We have
\bel{Theta diff est < flat}
\left|\left(\mathfrak{D}^\flat\right)^N \Theta^\nu_{\ell j}(x,y,\eta)\right|~\leq~\C_N~2^{-j\left({n-1\over 2}\right)},\qquad N\ge1.
\eeq
Moreover,   by (\ref{intersection of cone norm<}), the support of $\Theta^\nu_{\ell j}(x,y,\eta)$  in $\eta$ has a volume bounded by $\C~2^j2^{(j+\ell)\left({n-1}\right)}$.

Let $\Omega^\nu_{\ell j}(x,y)$ be given in (\ref{Omega rewrite}).  By using (\ref{Theta diff est < flat}), an $N$-fold integration by parts associated to $\mathfrak{D}^\flat$ shows that
\bel{Omega rewrite norm < flat}
\begin{array}{lr}\ds
\left|\Omega^\nu_{\ell j}(x,y)\right|~\leq~\C_{N}~2^{-j\left({n-1\over 2}\right)} ~2^j2^{(j+\ell)\left({n-1}\right)}
\\\\ \ds
\left\{1+4\pi^2 2^{2j}\left(\nabla_\eta\Phi\left(x,\L_\nu\eta_j^\nu\right)-\L_\nu^T y\right)_1^2+4\pi^2 2^{2(j+\ell)}\sum_{i=2}^{n}\left(\nabla_\eta\Phi\left(x,\L_\nu\eta_j^\nu\right)-\L_\nu^T y\right)_i^2 \right\}^{-N}. 
\end{array}
\eeq

Denote $\mathcal{X}=\mathcal{X}(x)=\left(\L_\nu^T\right)^{-1}\nabla_\eta\Phi\left(x,\L_\nu\eta_j^\nu\right)$.
By using (\ref{Omega rewrite norm < flat}), we have
\bel{int Omega rewrite norm < flat}
\begin{array}{lr}\ds
\int_{\R^n}\left|\Omega^\nu_{\ell j}(x,y)\right|dx
\\\\ \ds
~\leq~\C_{\Phi~N}\int_{\R^n} 2^{-j\left({n-1\over 2}\right)}2^j2^{(j+\ell)\left({n-1}\right)}\left\{1+ 2^{2j}(\mathcal{X}-y)_1^2+2^{2(j+\ell)}\sum_{i=2}^{n}(\mathcal{X}-y)_i^2\right\}^{-N}d\mathcal{X} 
\\\\ \ds
~\leq~\C_{\Phi~N}\int_{\R^n} 2^{-j\left({n-1\over 2}\right)}\left\{1+ \mathcal{Z}_1^2+     \sum_{i=2}^{n} \mathcal{Z}_i^2 \right\}^{-N} d\mathcal{Z}
\\\\ \ds~~~~~~~~~~~~~~~~~~~~~~~~~~~~~~~~~~~~~
\hbox{\small{$ \mathcal{Z}_1=2^{j}(\mathcal{X}-y)_1$~~and~~ $\mathcal{Z}_i=2^{j+\ell}(\mathcal{X}-y)_i,~~i=2,\ldots,n$}} 
\\\\ \ds
~\leq~\C_{\sigma~\Phi}~2^{-j\left({n-1\over 2}\right)}, \qquad \hbox{\small{$N\ge{n+1\over 2}$.}}
\end{array}
\eeq

Recall from section 6. For $\ell<-j/2-3$, $\Lambda_\ell$ intersects  only $\Gamma^\nu_j$ whose central direction  is $\xi^\nu_j=(\pm1,0)\in\R\times\R^{n-1}$.

We thus have
\bel{Omega_lj Sum norm < flat}
\begin{array}{lr}\ds
\int_{\R^n}\left|\Omega_{\ell j}(x,y)\right|dx~\leq~\sum_{\nu~\colon~\xi^\nu_j=(\pm1,0)} \int_{\R^n}\left|\Omega^\nu_{\ell j}(x,y)\right|dx\qquad\hbox{\small{by (\ref{Omega_lj Sum})}}
\\\\ \ds ~~~~~~~~~~~~~~~~~~~~~~~~~~~~
~\leq~\C_{\sigma~\Phi}~ 2^{-j\left({n-1\over 2}\right)}
\qquad \hbox{\small{by (\ref{int Omega rewrite norm < flat})}}

\\\\ \ds ~~~~~~~~~~~~~~~~~~~~~~~~~~~~
~\leq~\C_{\sigma~\Phi}~2^{\left({n-1\over 2}\right)\ell}.\qquad (-j-3\leq\ell)
\end{array}
\eeq

Observe that every $\p_y$ acting on $\Omega^\nu_{\ell j}(x,y)$ defined in (\ref{Omega^v_lj}) or (\ref{Omega rewrite}) gains a factor of $\C 2^j$ whenever $2^{j-1}\leq|\xi|=|\eta|\leq2^{j+1}$. 
By carrying out the same estimates in (\ref{Diff Ineq sigma eta_imath <})-(\ref{Omega_lj Sum norm < flat}), we find
\bel{nabla Omega_lj<} 
\int_{\R^n} \left|\nabla_y\Omega_{\ell j}(x,y)\right| dx~\leq~\C_{\sigma~\Phi}~
 2^j2^{\left({n-1\over 2}\right)\ell}
  \eeq
which further implies 
\bel{Est2 <} 
\int_{\R^n} \left|\Omega_{\ell j}(x,y)-\Omega_{\ell j}(x,x_o)\right| dx~\leq~\C_{\sigma~\Phi}~2^j|y-x_o|
2^{\left({n-1\over 2}\right)\ell}.
\eeq
\v

Recall $\mathfrak{Q}_r(x_o)$ defined in (\ref{rectangle R<})-(\ref{Q_r<}). Let $2^{k}\leq r^{-1}\leq 2^{k+1}$.
For $x\in{^c}\mathfrak{Q}_r(x_o)$, we must have
\bel{x-Phi Est <}
\left|\left(\L_\nu^Tx_o-\nabla_\eta\Phi\left(x,\L_\nu \eta^\nu_j\right)\right)_1\right|~\ge~2\cdot2^{-k}.
\eeq
If $y\in B_r(x_o)$, then $|y-x_o|\leq 2^{-k}$. For every $2^j\ge r^{-1}$,  we have
\bel{x,y-Phi Est <}
\begin{array}{lr}\ds
2^{2j}\left(\L_\nu^Ty-\nabla_\eta\Phi\left(x,\L_\nu \eta^\nu_j\right)\right)_1^2~\ge~ 2^{2(j-k)}
\\\\ \ds~~~~~~~~~~~~~~~~~~~~~~~~~~~~~~~~~~~~~~~~~~~~
~\ge~2^{j-k}.\qquad (j-k\ge0)
\end{array}
\eeq
By carrying out the same estimates in (\ref{Diff Ineq sigma eta_imath <})-(\ref{Omega_lj Sum norm < flat}), except for (\ref{Omega rewrite norm < sharp}) and (\ref{Omega rewrite norm < flat}) replaced  with
\bel{Omega rewrite norm < sharp*}
\begin{array}{lr}\ds
\left|\Omega^\nu_{\ell j}(x,y)\right|
~\leq~\C_{N}~2^{-j\left({n-1\over 2}\right)}~2^j2^{(j+\ell)\left({n-1}\right)}
~2^{-j+k}
\\\\ \ds~~~~~~~~~~~~~~~~~~~~~~~~
\Bigg\{1+4\pi^2 2^{2j}\left(\nabla_\eta\Phi\left(x,\L_\nu\eta_j^\nu\right)-\L_\nu^T y\right)_1^2 \Bigg\}^{1-N} 
~~~~~~~
\hbox{\small{by (\ref{x,y-Phi Est <})}}
\end{array}
\eeq
 and
\bel{Omega rewrite norm < flat*}
\begin{array}{lr}\ds
\left|\Omega^\nu_{\ell j}(x,y)\right|~\leq~\C_{N} ~2^{-j\left({n-1\over 2}\right)}~~2^j2^{(j+\ell)\left({n-1}\right)}
~2^{-j+k}
\\\\ \ds
\left\{1+4\pi^2 2^{2j}\left(\nabla_\eta\Phi\left(x,\L_\nu\eta_j^\nu\right)-\L_\nu^T y\right)_1^2+4\pi^2 2^{2(j+\ell)}\sum_{i=2}^{n}\left(\nabla_\eta\Phi\left(x,\L_\nu\eta_j^\nu\right)-\L_\nu^T y\right)_i^2 \right\}^{1-N}
\\ \ds
\qquad~~~~~~~~~~~~~~~~~~~~~~~~~~~~~~~~~~~~~~~~~~~~~~~~~~~~~~~~~~~~~~~~~~~~~~~~~~~~~~~~~~~~~~~~~~~~\hbox{\small{by (\ref{x,y-Phi Est <})}}
\end{array}
\eeq
we find
\bel{Est3 <} 
\int_{{^c}\mathfrak{Q}_r(x_o)} \left|\Omega_{\ell j}(x,y)\right| dx~\leq~\C_{\sigma~\Phi}~ {2^{-j}\over r}
~ 2^{\left({n-1\over 2}\right)\ell}
 \eeq
for every $y\in B_r(x_o)$ whenever $2^{j}\ge r^{-1}$.

Up to this end, we have obtained the desired estimates  in (\ref{Est1*})-(\ref{Est3*}).

On the other hand, consider $\Omega^*_{\ell j}(x,y)$  defined in (\ref{Omega^*_t,j}) associated to the adjoint operator $\F^*_\ell$   in (\ref{F*_t singular}).   
We prove  $\Omega^*_{\ell j}(x,y)$  satisfying (\ref{Est1*}), (\ref{Est2*}) and (\ref{Est3*}) where $\Q_r(x_o)$ is  replaced by $\Q^*_r(x_o)$ defined  in (\ref{rectangle R*})-(\ref{Q_r^*<}), by carrying out the same estimates  in  {\bf Case One}, {\bf Two} and {\bf Three}, with some necessary changes. In particular,   we shall apply  (\ref{rectangle R*}), (\ref{rectangle R*>}),  (\ref{rectangle R*<}) instead of (\ref{rectangle R}), (\ref{rectangle R>}), (\ref{rectangle R<}) respectively.

Lastly, we make a final remark on a possible extension of our result to the $n$-parameter setting. 
We say $\sigma\in\Hat{\S}^m$ if 
\bel{Class n-parameter}
\left|\p_\xi^\alpha\p_{x,y}^\beta \sigma(x,y,\xi)\right|~\leq~\C_{\alpha~\beta}~\left(1+|\xi|\right)^m\prod_{i=1}^n \left({1\over 1+|\xi_i|}\right)^{\alpha_i}
\eeq
for every multi-indices $\alpha,\beta$.
\v
{\bf Conjecture A*} ~{\it Let $\sigma\in \Hat{\S}^m$  as (\ref{Class n-parameter}) for $-(n-1)/2<m\leq0$. Fourier integral operator $\F$   defined in (\ref{Ff})-(\ref{nondegeneracy}) extends to a bounded operator 
  \bel{Result n}\left\| \F f\right\|_{\L^p(\R^n)}~\leq~\C_{p~\sigma~\Phi}~\left\| f\right\|_{\L^p(\R^n)}
 \eeq
whenever 
\bel{Formula n}
\left| {1\over 2}-{1\over p}\right|~\leq~{-m\over n-1}.
\eeq}
\v
Certain preliminary results have been developed in the same spirt of cone decomposition. 
Essentially, in this $n$-parameter setting, it can be understood as an $n$-fold covering by the $2$-parameter cone decomposition, introduced in section 3,  which already carries the main characteristic of our framework. For this reason, we choose to focus on the $2$-parameter Fourier integral operators in the present paper.

\v

{\bf Acknowledgement:

~~~~~~~~~~~~~~~~~~~~~~~~~I am deeply grateful to my advisor Elias M. Stein for those unforgettable lectures given in Princeton.}

email: wangzipeng@westlake.edu.cn

\end{document}